\documentclass[12pt]{article}
\usepackage{amstext,amsfonts,amssymb,amsmath,amsopn}
\usepackage[frenchb,english]{babel}
\usepackage{color}

\newcommand{\Z}{\mbox{${\mathbb Z}$}}
\newcommand{\Q}{\mbox{${\mathbb Q}$}}

\newcommand{\F}{\mbox{${\mathcal F}$}}
\newcommand{\FF}{\mbox{${\mathbb F}$}}

\newcommand{\lito}{\mbox{o}}
\newcommand{\bul}{\textbullet\ }

\newcounter{tab}
\def\tabcaption #1
  {
   \begin{center}
   {\sc Table} \arabic{tab} {\it #1} 
    \end{center}\addtocounter{tab}{1}
}

\newtheorem{prop}{Proposition}
\newtheorem{est}{Estimation num\'erique}

\begin{document}

\selectlanguage{frenchb}

\bigskip

\begin{center}
{\Large\bf Polyn\^omes de degr\'e sup\'erieur \`a 2
prenant beaucoup de valeurs premi\`eres}

Fran\c cois DRESS et Bernard LANDREAU
\footnote{Avec le concours du 
Centre de Calcul Intensif des Pays de Loire
et du Groupement De Services Mathrice du CNRS}

\bigskip
\centerline{\it 7 juin 2012}
\bigskip

\end{center}

\selectlanguage{english}
\begin{abstract}

	For degrees $3$ to $6$, we first give numerical results on polynomials which take many prime values on an interval of consecutive values of the variable.

	In particular, we have improved Ruby's record for the "$n$
        out of $n$" case, for $n\ =\ 58$, by using a polynomial of degree $6$.

	In the theoretical part of this paper, we describe a
        heuristic probabilistic model in the "$n$ out of $n$" case:
        exactly $n$ (different) prime values on an interval of $n$ consecutive values of the variable. We find that the heuristic value of the probability of the event "$n$ out of $n$" for a generic polynomial is equal to the product of two factors: an arithmetic factor related to global conditions of non-divisibility, and a size factor determined by the position of the polynomial in a "well-shaped" domain of the space of coefficients. This leads to a heuristic estimate for the number of "$n$ out of $n$" polynomials in a given "well-shaped" domain. Finally, results of extended numerical experiments show a satisfactory agreement with the heuristic values given by the model.
      \end{abstract}

\selectlanguage{frenchb}

\section{Introduction}

L'histoire des polyn\^omes qui prennent "beaucoup" de valeurs premi\`eres commence au polyn\^ome $x^2 + x + 41$ d'Euler (1772), qui prend 40 valeurs premi\`eres (distinctes) de $x = 0$ \`a $x = 39$. 
Le crit\`ere de Rabinowitch (\cite{Rab}, 1912) fournit une raison alg\'ebrique
\`a cette performance, et en m\^eme temps il \'enonce que, pour les polyn\^omes de la forme $x^2 + x + C$,
41 est la derni\`ere valeur op\'erante de C.

L'arriv\'ee des ordinateurs a modifi\'e la situation en entra\^{\i}nant dans les ann\'ees 1980 des recherches exp\'erimentales. Deux types de probl\`emes ont alors \'et\'e clairement distingu\'es : 
\begin{itemize}
\item[-] le probl\`eme $n$ sur $n$ : il s'agit de trouver des polyn\^omes qui prennent $n$ valeurs premi\`eres sur une s\'equence de $n$ valeurs cons\'ecutives de la variable,
\item[-] le probl\`eme $k$ sur $n$ : il s'agit de trouver des polyn\^omes qui prennent un grand nombre $k$ de valeurs premi\`eres sur une s\'equence de $n$ valeurs cons\'ecutives de la variable. 
\end{itemize}

Des records du type "$k$ sur 1\,000" ont d'abord \'et\'e obtenus --- voir Ribenboim (\cite{Rib}, 1996).
 Ensuite Fung (1988, cit\'e par \cite{Mol}) suivi par Ruby (1989, cit\'e \'egalement par \cite{Mol}) ont trouv\'e des polyn\^omes donnant 43 puis 45 valeurs premi\`eres (distinctes) pour autant de valeurs cons\'ecutives  de la variable. 
Ces records ayant sembl\'e imbattables (et ils le sont tr\`es
vraisemblablement si l'on se limite au degr\'e 2), les recherches
exp\'erimentales se sont alors concentr\'ees sur ce que nous appellons
le cas "$k$ sur $n$", 
en degr\'e 2 \`a l'unique exception d'un r\'esultat de Goetgheluck
(\cite{Goe},1989) en degr\'e 3.
La situation est rest\'ee confuse quelques ann\'ees jusqu'\`a
l'article de Boston et Greenwood ([2], 1995) qui normalise de fa\c con
efficace les polyn\^omes (et donne une liste de "bons" polyn\^omes
pour $n = 100$).

Un article de Dress et Olivier ([3], 1999) explore tr\`es
compl\`etement le cas du degr\'e 2 : il donne un mod\`ele probabiliste
heuristique et met en \'evidence ce qu'ils appellent le mur de Schinzel, il montre la diff\'erence 
entre la probl\'ematique des records "$n$ sur $n$" et celle des records "$k$ sur $n$", 
et enfin il donne des r\'esultats num\'eriques tr\`es \'etendus.

L'objet de notre article est de prolonger l'\'etude aux degr\'es sup\'erieurs (de 3 \`a 6) :
 exploration num\'erique g\'en\'erale, et extension du mod\`ele probabiliste heuristique dans le cas "$n$ sur $n$". 
 Un nouveau record "$n$ sur $n$" est \'etabli : le polyn\^ome
$$\frac{1}{72}x^6+\frac1{24}x^5-\frac{1\,583}{72}x^4-\frac{3\,161}{24}x^3+\frac{200\,807}{36}x^2+\frac{97\,973}3x-11\,351$$
prend 58 valeurs premi\`eres pour 58 valeurs cons\'ecutives de la
variable, de $x = -25$ \`a $x = 31$.

\section{Probl\'ematique g\'en\'erale} 

\subsection{Normalisation des polyn\^omes} 

   Il faut tout d'abord donner les conventions qui pr\'ecisent l'objet des recherches :
\begin{itemize}
\item [-] les valeurs premi\`eres compt\'ees doivent avoir des valeurs absolues toutes distinctes,
\item [-] 1 est accept\'e comme nombre premier exceptionnel (cela arrive tr\`es rarement, et  la pr\'esence d'une valeur 1 ou $-1$ sera chaque fois explicitement signal\'ee), 
\item [-] les polyn\^omes  sont \`a coefficients rationnels, pas
  n\'ecessairement unitaires et prennent des valeurs enti\`eres sur \Z.
\end{itemize}

   Les notations et la normalisation des polyn\^omes sont
destin\'ees \`a \'eviter les redondances dues aux
transformations du type $\pm P(\pm x - m)$. 
\begin{itemize}
\item[-] Pour les polyn\^omes \`a coefficients entiers, on suit la
normalisation de Boston et Greenwood : $P$ s'\'ecrit
$$
P(x) = a_0x^d + a_1x^{d-1} +\cdots + a_{d-1}x + a_d ,\qquad
\mbox{ avec } a_0 \ge 1 \quad \mbox{ et }\quad  0 \le a_1 \le \frac{d}2 a_0 ,
$$
\item[-] Pour les polyn\^omes \`a coefficients rationnels, on peut \'ecrire les polyn\^omes dans la
base naturelle des polyn\^omes \`a valeurs enti\`eres, 
$$ P(x) = b_0\frac{x(x-1)\ldots(x-d+1)}{d!} +
b_1\frac{x(x-1)\ldots(x-d+2)}{(d-1)!} +\cdots + b_{d-1}x + b_d ,$$
avec les $b_j$ entiers. Mais il faut alors les r\'e\'ecrire
sous la forme $c_0 x^d + c_1 x^{d-1} + \cdots + c_{d-1} x + c_d$ avec
les $c_j$ rationnels de d\'enominateur $d!$
pour pouvoir appliquer 
la normalisation de Boston et Greenwood.
Cette  normalisation conduit alors \`a imposer $c_0 \ge 1$ et $ 0 \le c_1 \le \frac{d}2 c_0$.
\end{itemize}

\subsection{Les diviseurs premiers p\'eriodiques} 

   Pour simplifier l'exposition des r\'esultats et des  conjectures, 
Dress et Olivier ont introduit la notion de {\em diviseur premier
p\'eriodique}, en abr\'eg\'e d.p.p., d'un polyn\^ome $P$ : 
on d\'esigne ainsi tout nombre premier $p$ qui, pour tout $m$, 
divise au moins l'un des entiers $P(m), P(m+1), ..., P(m+p-1)$
(cela arrive d\`es qu'il existe $m \in\Z$ tel que $p | P(m)$, 
tout diviseur premier est automatiquement p\'eriodique).

   Pour les performances "$n$ sur $n$", l'absence des d.p.p. inf\'erieurs 
ou \'egaux \`a $n$ est obligatoire. 
Pour les records "$k$ sur $n$", leur influence est plus modul\'ee,
les d.p.p. seront les nombres premiers qui fournissent des facteurs
non triviaux dans le produit infini qui donne la constante 
de Hardy-Littlewood de $P$, voir \cite{Har}.

\subsection{Polyn\^omes \`a coefficients entiers ou rationnels ?} 

   Les polyn\^omes \`a coefficients rationnels et \`a valeurs enti\`eres
sur \Z\  sont les combinaisons lin\'eaires \`a coefficients entiers des polyn\^omes
$$
1, x,\frac{x(x-1)}{2} ,\ldots ,\frac{x(x-1)\ldots (x-k+1)}{k!}
,\ldots
$$

   En degr\'e 2 ou 3, tout polyn\^ome \`a coefficients rationnels
qui n'est pas dans $\Z[x]$ admet forc\'ement 2 ou 3 comme d.p.p.,
ce qui le rend tout \`a fait inint\'eressant pour notre recherche.

   En revanche, \`a partir du degr\'e 4, on voit appara\^{\i}tre des
polyn\^omes \`a coefficients rationnels qui ne sont pas dans $\Z[x]$
et qui peuvent \^etre de bons candidats. 
Cela accro\^{\i}t consid\'erablement le "vivier" de polyn\^omes \`a
tester. Une recherche a donc \'et\'e effectu\'ee sur ces polyn\^omes,
cela a permis de battre nettement les records obtenus sur les
polyn\^omes \`a coefficients entiers.

\subsection{La conjecture de Schinzel} 

Si $P$ est un polynome (\`a valeurs enti\`eres), on notera par $np(P)$
le maximum du nombre de valeurs cons\'ecutives de $P$ toutes
premi\`eres (et distinctes en valeur absolue).
 On a alors toujours $np(P)\le
2p-1$ o\`u $p$ est le plus petit d.p.p de $P$. 
Si on fixe un entier $d$, on notera $np_d$ le supremum de $np(P)$ sur tous les polyn\^omes $P$ de degre $d$.

La conjecture de Schinzel ([12]) peut s'\'enoncer comme suit  dans le cas d'un unique polyn\^ome $P$ (\`a coefficients entiers) : si $m$ est inf\'erieur au  plus petit d.p.p. de $P$, alors il existe une infinit\'e de $m$ tels que les entiers $P(m), P(m+1),\ldots, P(m+n-1)$ soient tous premiers.
	
Soient $P$ un polyn\^ome et $p$ son plus petit d.p.p., la conjecture
de Schinzel implique qu'en g\'en\'eral  $np(P) =p-1$ 
 (en n\'egligeant le cas, 
 rare en pratique,  o\`u $p$ ou $-p$ serait une des valeurs
 prises dans une suite de longueur maximale de valeurs cons\'ecutives
 toutes premi\`eres).

Nous avons effectu\'e dans \cite{Dre} une construction par
congruences sur les  coefficients d'un polyn\^ome, avec utilisation du
th\'er\`eme chinois,  qui permettait d'obtenir  un polyn\^ome n'ayant
aucun d.p.p. inf\'erieur ou \'egal \`a une valeur $z$ fix\'ee. 
On en d\'eduit que la conjecture de Schinzel implique que $np_d = + \infty$ pour tout $d$.

Un r\'esultat tres r\'ecent de Granville (\cite{Gra}) d\'ecoulant du th\'eor\`eme
 de Green et Tao (\cite{Gre}) sur les progressions arithm\'etiques dans la
 suite des nombres premiers assure que, pour tout $d\ge1$ et tout $n$,
 il existe une infinit\'e de polyn\^omes $P$ tels que $P(0),
 P(1),\ldots, P(n-1)$ soient premiers distincts.
 Cela entra\^{\i}ne
notamment que $np_d = +\infty$ pour tout $d$ de fa\c con s\^ure.
De plus, Granville conjecture l'existence pour tout $d\ge1$ et tout
$n$, d'un polyn\^ome unitaire $P$ v\'erifiant $0 < P(0) < P(1) <
\cdots < P(n-1)$ tous premiers et $P(n-1)$ de l'ordre de
$\left((1 + \lito(1))\frac{n}{e^\gamma d}\right)^{n/d}$ pour $n>=4(d \log d)^2$.

\subsection{Le mur de Schinzel et la zone observable} 

L'\'etude probabiliste heuristique men\'ee par Dress et
Olivier en [3]
 met en lumi\`ere un ph\'enom\`ene baptis\'e  "mur de
Schinzel" qui renvoie le fonctionnement effectif de la conjecture de Schinzel  \`a des zones de
valeurs qui sont hors de toute possibilit\'e d'exp\'erimentation
num\'erique. La taille heuristique du plus petit entier $m$
mentionn\'e dans la conjecture de Schinzel est en effet en g\'en\'eral
\'enorme,  de l'ordre en degr\'e $d$ de $n^{\frac{d}2n}$ ,
 et tr\`es vite hors d'atteinte de tout calcul explicite.
La borne impliqu\'ee par le r\'esultat de Granville, quoique tr\`es
inf\'erieure \`a la valeur du mur de Schinzel, reste pareillement hors
d'atteinte hors d'atteinte exp\'erimentale (en outre la condition sur $n$ est un peu forte, 
par exemple elle s'\'ecrit $n\ge 133$ pour $d = 4$, mais ce point est
tout \`a fait mineur).
\\
On constate exp\'erimentalement que pour un polyn\^ome donn\'e $P(x)$
 sans obstructions arithm\'etiques dues \`a des petits d.p.p.,
et avec des coefficients ``raisonnablement'' petits, la meilleure
performance $n$ sur $n$ en fonction de la variable est atteinte avant
$x=2n$ en valeur absolue. Le mod\`ele heuristique confirme cette
observation : si aucune performance n'a \'et\'e constat\'ee au
voisinage de l'origine, alors il faudra attendre le mur de Schinzel
(ou peut-\^etre la borne de Granville)
pour voir une performance avec la valeur conjecturale maximale de $n$
se produire. 
Cela justifie donc, malgr\'e la conjecture de Schinzel, la recherche
de polyn\^omes records sur une zone observable \`a l'\'echelle de nos
ordinateurs.

  On peut se poser la question d'une d\'efinition pr\'ecise de la zone
observable. Cela ne peut \^etre qu'arbitraire, et
sans v\'eritable enjeu, tellement l'\'ecart entre les valeurs
"accessibles" et le mur de Schinzel (et m\^eme la borne de Granville) 
devient gigantesque \`a partir de $n = 20$ environ. 
Nous proposons la valeur $n^{d\sqrt{n}}$ comme borne en
degr\'e $d$ \`a la fois pour les valeurs (absolues) prises par le
polyn\^ome, et pour le num\'erateur et le d\'enominateur des
coefficients (on peut noter que, pour le polyn\^ome d'interpolation de 
$n$ valeurs enti\`eres distinctes arbitraires, 
le maximum des num\'erateurs et d\'enominateurs des coefficients est
de l'ordre de $2^n$). 
On peut maintenant d\'efinir $np^*(P)$ comme le maximum du nombre de
valeurs cons\'ecutives toutes premi\`eres, sous les restrictions
\'enonc\'ees ci-dessus de borne $n^{d\sqrt{n}}$, puis $np_d^*$ comme le maximum 
de $np^*(P)$ sur tous les polyn\^omes de degr\'e $d$,
 et sous les m\^emes restrictions. 
Malgr\'e les conditions suppl\'ementaires impos\'ees, il est
impossible de calculer exactement la valeur de $np^*(P)$ pour un
polyn\^ome donn\'e, mais on peut exp\'erimentalement en proposer une 
valeur pr\'esum\'ee, extr\^emement probable. Ce sont ces valeurs
pr\'esum\'ees (en toute rigueur des minorations) que nous donnons ici,
et  les valeurs \'egalement pr\'esum\'ees
 de $np_d^*$ qui en r\'esultent.

   Enfin, on peut d\'efinir de fa\c con exactement similaire
 pour les nombres de valeurs premi\`eres 
sur $n$ valeurs cons\'ecutives, les maximums $nk(P,n)$ et
   $nk^*(P, n)$, $nk_d(n)$ et
$nk_d^*(n)$.  
\vspace{-3mm}

\section{R\'esultats exp\'erimentaux} 

\subsection{Le probl\`eme " $n$ sur $n$ "} 

  	Les statistiques effectu\'ees sur les observations exp\'erimentales, confirm\'ees par le mod\`ele probabiliste heuristique d\'ecrit au paragraphe 4, permettent d'identifier les facteurs de succ\`es de la recherche num\'erique, qui reposent respectivement sur l'arithm\'etique, la taille, et la "chance".

	Notons tout d'abord que, pour qu'un polyn\^ome puisse donner une performance 
"$n$ sur $n$", il faut que son plus petit d.p.p. soit sup\'erieur \`a
$n$, sauf dans le cas particulier et rare en pratique o\`u $P$
prendrait la valeur $p$ ou $-p$.

	Le facteur arithm\'etique est la restriction (par un criblage efficace) aux polyn\^omes n'ayant pas de petits d.p.p., et on \'enoncera dans le paragraphe 4 sur la mod\'elisation un th\'eor\`eme arithm\'etique qui montre que le gain apport\'e par cette restriction est beaucoup plus important qu'on ne pouvait attendre a priori.

	Le facteur taille est la limitation de la recherche aux plages
        de la variable proches de l'origine avec des polyn\^omes \`a
        petits coefficients. Cela permet de maximiser la valeur  qui
        repr\'esente (hors facteur arithm\'etique) la probabilit\'e
        heuristique que $P(m)$ soit premier. C'est probablement ce
        deuxi\`eme facteur qui explique pour le probl\`eme "$n$ sur
        $n$" l'am\'elioration des performances lorsque l'on \'el\`eve
        le degr\'e (les meilleurs polyn\^omes oscillent et restent
        plus longtemps avec des "petites" valeurs).

	Nous avons appel\'e "chance" le troisi\`eme facteur. En effet, on constate dans les recherches exp\'erimentales \'etendues, qu'une fois les deux premiers facteurs pris en compte, la complexit\'e pourtant toute d\'eterministe du probl\`eme donne l'apparence du hasard. Le seul moyen d'action sur ce facteur est alors le nombre de polyn\^omes test\'es.

	NOTA sur les records  et leurs auteurs.

Except\'e ceux d'Euler, Fung et Ruby, les records list\'es ci-dessous
jusqu'au degr\'e 6 ont \'et\'e d\'ecouverts par Dress et Landreau dans
le cadre du pr\'esent travail. Trois records principaux (46 en
degr\'e 3, 46 en degr\'e 4 \`a coefficients tous entiers, et 57 en
degr\'e 5 \`a coefficients rationnels non tous entiers) ont \'et\'e
pr\'e-publi\'es dans \cite{Rib} :
 P. Ribenboim, The little book of bigger primes, 2d edition, Springer,
 2004, p. 148.
   Par
 contre, en degr\'e 6, le tout dernier record 58 d\'ecouvert par Dress et
 Landreau (2010) est publi\'e pour la premi\`ere fois ici.
\\
Certains polyn\^omes ont \'et\'e red\'ecouverts dans la comp\'etition Internet \cite{AZPC} : Al
 Zimmermann's Programming Contest, Prime Generating Polynomials,
 organis\'ee par Ed Pegg Jr en juillet  2006. Cela sera mentionn\'e.
	Signalons enfin deux particularit\'es de cette comp\'etition. Primo une cat\'egorie sp\'eciale \'etait consacr\'ee aux polyn\^omes prenant (sur la zone record) des valeurs toutes de m\^eme signe, qui ne seront pas donn\'es ici. Secundo il n'\'etait pas impos\'e que les valeurs premi\`eres (obtenues pour des valeurs cons\'ecutives de la variable) soient toutes diff\'erentes en valeur absolue ; cette sp\'ecification \'etait impos\'ee par Boston et Greenwood, que nous avons suivi pour les records (mais non pour l'heuristique). 

	{\bf Degr\'e 2}
\\[1ex]
$n = 45$ :   1 polyn\^ome :   $36x^2 + 18x -1\, 801$ (polyn\^ome de
Ruby), de $x = -33 $ \`a $x = 11$
\\[1ex]
$n = 43$ :   2 polyn\^omes : $47x^2 + 9x - 5 209$ (polyn\^ome de Fung), de $x = -22$ \`a $x = 18$\\
	\hbox{}	\hspace{4cm}				      $103 x^2 + 31x  - 3\, 391$ (Ruby)
\\[1ex]
$n = 40$ :   7 polyn\^omes   (dont le polyn\^ome d'Euler)

	{\bf Degr\'e 3}
\\[1ex]
$n = 46$ :   1 polyn\^ome :   $6x^3 + 83x^2 - 13\,735x + 30\, 139$, de $x=-26$ \`a $x = 19$ (DL)\\
$n = 41$ :   3 polyn\^omes\\
$n = 40$ :   7 polyn\^omes

	{\bf Degr\'e 4}
\\[1ex]
$n = 49$ :   2 polyn\^omes : $\frac34x^4 + \frac12 x^3 -\frac{4\,323}{4}  x^2 + \frac{34\,415}2x - 62\,099$
(DL, red\'ecouvert par \\ \hbox{}	\hspace{4cm} J. Wroblewski et J.-C. Meyrignac dans [AZPC])\\
	\hbox{}	\hspace{4cm} $\frac94x^4 +\frac52  x^3 - \frac{5077}4 x^2 - \frac{24\,951}2x - 347$                               (DL)\\
$n = 47$ :   1 polyn\^ome\\
$n = 46$ :   8 polyn\^omes \`a coefficients entiers  (DL)
		 + 1 polyn\^ome \`a coefficients rationnels non entiers

	{\bf Degr\'e 5}
\\[1ex]
$n= 57$ : 1 polyn\^ome : 
$\frac14 x^5 +\frac12 x^4 -\frac{345}4 x^3 +\frac{879}2 x^2 +  17\, 500 x + 70\, 123$  
(DL, red\'ecouvert \\ \hbox{}	\hspace{4cm} par Shyam Sunder Gupta dans [AZPC])
\\
$n = 51$ : 1 polyn\^ome : $\frac12 x^5 + \frac34 x^4 - 49 x^3
-\frac{463}4 x^2 +\frac{15\,099}2  x + 3\, 457$  (DL)\\
$n = 50$ :   4 polyn\^omes \`a coefficients rationnels non entiers (dont un prend la valeur -1)\\
$n = 49$ :   1 polyn\^ome \`a coefficients entiers :
 $3 x^5  + 7 x^4  - 340 x^3  - 122 x^2  + 3876 x + 997$ (DL 2001),
	 + plusieurs polyn\^omes \`a coefficients rationnels non entiers

	{\bf Degr\'e 6}
\\[1ex]
$n = 58$ :   1 polyn\^ome  :
$\frac{1}{72}x^6+\frac1{24}x^5-\frac{1\,583}{72}x^4-\frac{3\,161}{24}x^3+\frac{200\,807}{36}x^2+\frac{97\,973}3x-11\,351$ (DL 2010)  
\\[1ex]
$n = 57$ :   2 polyn\^omes  :
$\frac1{36}x^6+\frac1{12}x^5-\frac{125}9x^4-\frac{3\,791}{12}x^3-\frac{76\,829}{36}x^2-\frac{15\,277}6x-58\,567$
\\[1ex]
et
$\frac1{72}x^6+\frac1{24}x^5-\frac{1343}{72}x^4-\frac{1265}{24}x^3+\frac{158495}{36}x^2+\frac{6044}3x-113\,
723$
(DL 2009 et 2010)
\\[1ex]
$n = 55$ :   1 polyn\^ome :  
$\frac1{36}x^6 -\frac{199}{18}  x^4 -\frac{71}2  x^3 +\frac{43\,165}{36}  x^2 +\frac{11\,639}2 x -  2\,423$ 
\\	\hbox{}	\hspace{5cm} (J. Wroblewski et J.-C. Meyrignac dans [AZPC])
\\
$n = 54$ :   3 polyn\^omes :  
\\
$n = 53$ :   10 polyn\^omes   (DL)
\\
$n = 44$ :   2 polyn\^omes \`a coefficients entiers :\\
$x^6 - 107 x^5 + 4\,967x^4 - 108\,362x^3 + 1\,387\,098 x^2 - 9\,351\,881
x + 25\,975\,867$ (non normalis\'e)
\\	\hbox{}	\hspace{5cm}(J. Wroblewski et J.-C. Meyrignac dans [AZPC])
\\
$x^6+2x^5-100x^4+79x^3+367x^2-3\, 919x-4\, 723$ (DL 2011)

\vspace{-5mm}
$$
\begin{array}{|c|c|c|c|c|c|}
\hline
d & 2 & 3 & 4 & 5 & 6 \\
\hline
\mbox{minoration (coefficients tous entiers)} & 45 & 46 & 46 & 49 &
\it{44} \\
\hline
\mbox{minoration (coefficients rationnels non tous entiers)} & & & 49
& 57 & \it{58}\\
\hline
\end{array}
$$
\tabcaption{minorations exp\'erimentales de $np^*_d$\\ (qui sont en m\^eme temps les valeurs pr\'esum\'ees si $d\le 5$)}

\subsection{Le probl\`eme "$k$ sur $n$"} 

  	Les facteurs de succ\`es de la recherche num\'erique sont bien s\^ur les m\^emes que ceux du probl\`eme "$n$ sur $n$", mais leur mode d'action est tr\`es diff\'erent. En cons\'equence, la mod\'elisation heuristique est \'egalement tr\`es diff\'erente (mais nous ne l'exposerons pas dans cet article).

	Le facteur arithm\'etique ne se limite pas \`a l'\'elimination
        des polyn\^omes \`a petits d.p.p. Pour le probl\`eme "$k$ sur
        $n$", il faut consid\'erer globalement un grand nombre des
        "premiers" d.p.p. du polyn\^ome, dont la r\'epartition va
        fortement influencer ses performances. 
Une conjecture due \`a Hardy et Littlewood (\cite{Har}, 1923)
\'enonce que le nombre $\pi_{P}(x)$ des valeurs premi\`eres prises
jusqu'\`a $x$ par un polyn\^ome $P$ est \'egal au produit du logarithme int\'egral
de $x$ par la ``constante de Hardy-Littlewood'', qui est un produit
infini dont le facteur g\'en\'erique essentiel est $( 1-\frac{\omega(p)}{p})$,
o\`u $\omega(p)$ d\'esigne le nombre de solutions de la congruence $P(x) \equiv 0 \pmod{p}$.

	Cette conjecture de Hardy et Littlewood \'etait limit\'ee au
        degr\'e 2 mais elle peut \^etre \'etendue aux degr\'es
        sup\'erieurs (Bateman et Horn, \cite{Bat}, 1962).
 Enfin, dans le cas de polyn\^omes \`a coefficients rationnels, 
elle peut \'egalement s'\'etendre en rempla\c cant le quotient $\omega(p)/p$ par le
quotient $\omega(p,t)/p^t$, o\`u $p^t$ ($t$ entier $\ge1$) est la
p\'eriode de la suite $P(n)$ modulo $p$ et 
 $\omega(p,t)$  le  nombre de solutions de la congruence $P(x)\equiv 0 \pmod{p^t}$.

Finalement, les performances d'un polyn\^ome donn\'e ont heuristiquement et exp\'erimentalement
un comportement de loi binomiale avec une probabilit\'e li\'ee \`a sa constante de Hardy-
Littlewood. Les meilleures performances en "$k$ sur $n$" sont obtenues pour des polyn\^omes
ayant une constante de Hardy-Littlewood $C(P)$ \'elev\'ee (voir par exemple Fung et
Williams, \cite{Fun}, 1962, et Jacobson et Williams, \cite{Jac}, 2003, dans une probl\'ematique plus
th\'eorique que la n\^otre).
 Nous avons \'etudi\'e les performances du polyn\^ome de Jacobson et
Williams record pour $C(P)$ en discriminant positif, $x^2 + x - A$,
avec $A = 1\,231\,847\,748\,861\,730\,729$
 ($C(P) = 5.24376$). La conjecture d'Hardy-Littlewood est trop globale mais elle peut
s'interpr\'eter localement par une probabilit\'e heuristique $p$ que $P(x)$ soit premier, voisine
de $C(P)/(2 \log |P(x)|)$. On en d\'eduit alors qu'heuristiquement les records sont \`a chercher
dans une zone tr\`es limit\'ee autour de $x = \sqrt{A} = 1\,109\,886\,368$. Les calculs heuristiques selon
la m\'ethode du paragraphe 4.5 du pr\'esent article sugg\`erent qu'on pourrait trouver le
record \`a 17 ou 18 : de fait, nous avons trouv\'e 17 valeurs cons\'ecutives premi\`eres \`a partir
de $x_0 = 1\,925\,947\,321$.
 On peut \'egalement chercher le record $k$ sur 1 000. L'heuristique
repose alors (comme dans \cite{Dre}) sur une loi binomiale et les calculs sugg\`erent qu'on
pourrait trouver le record vers 395 (valeur de $\mu n + 2\sigma \sqrt{n}$ pour $n = 1 000$) : de fait, nous
avons trouv\'e 399 valeurs premi\`eres sur 1 000 \`a partir de $x_0 = 1\,109\,886\,131$. Ces r\'esultats
num\'eriques sur le polyn\^ome de Jacobson et Williams appellent deux commentaires :
d'une part ils corroborent notre heuristique, d'autre part ils montrent qu'il ne suffit pas
d'une bonne constante de Hardy-Littlewood si on le "paie" d'une augmentation excessive
des valeurs prises (on retrouve le facteur "chance").

En ce qui nous concerne, la constante $C(P)$ n'est jamais calcul\'ee pr\'ecis\'ement dans la
recherche exp\'erimentale car les temps de calcul sont prohibitifs ; par contre la somme des
"premiers" $\omega(p)/p$ fournit un excellent crit\`ere empirique qui
est \`a la base du criblage que nous avons effectu\'e.

	Les recherches sont effectu\'ees sur des plages de la variable
        proches de l'origine, mais encore plus ici que pour le
        probl\`eme "$n$ sur $n$", le facteur heuristique li\'e \`a la
        taille des valeurs est crucial. Il s'ensuit que les meilleures
        performances sont obtenues, \`a degr\'e fix\'e, pour les
        polyn\^omes de petit coefficient directeur. Ce ph\'enom\`ene
        n'intervenait pas pour les performances "$n$ sur $n$", il est
        ici d'autant plus prononc\'e que $n$ est grand.
Cela \'etant, la d\'ependance pr\'ecise par rapport au degr\'e, qui
est un des aspects de ce facteur taille, est malais\'ee \`a saisir.

	Le facteur "chance", qui intervenait tr\`es brutalement dans le cas "$n$ sur $n$", est ici moins visible et fonctionne en deux temps : primo pour fournir apr\`es criblages des polyn\^omes ayant une grande constante de Hardy-Littlewood, secundo pour donner des performances avantageuses dans le comportement de loi binomiale \'evoqu\'e plus haut.

	Nous donnons ci-dessous, tout d'abord un tableau r\'ecapitulatif des performances (minorations exp\'erimentales de $nk_d^*(n)$, qui sont en m\^eme temps les valeurs pr\'esum\'ees sauf si $d$ et $n$ sont grands), ensuite la liste des polyn\^omes records.

 $$
       \begin{array}{|c|c|c|c|c|c|}
\hline
d \slash n & 50 & 100 & 200 & 500 & 1\, 000 \\
\hline
2 & 49 & 90 & 166 & 369 & 698 \\
\hline
3 & 48 & 87 & 154 & 338 & 601 \\
\hline
4 & 49 & 88 & 156 & 316  & 539\\
\hline
5 & 50 & \it{84} & \it{134} & \it{258} & \it{429} \\
\hline
6 & 50 & \it{83} & \it{132} & \it{261} & \it{404}\\
\hline
\end{array}
$$
\tabcaption{Meilleurs r\'esultats "$k$ sur $n$" en fonction du degr\'e}

\begin{center}
\it (les nombres \'ecrits en italiques concernent les zones o\`u les
recherches n'ont pas pu \^etre suffisamment pouss\'ees et o\`u les
r\'esultats ne sont donc vraisemblablement pas optimaux
- et nous ne donnerons d'ailleurs pas nos polyn\^omes records)
\end{center}

{\bf Degr\'e 2}
\\[1ex]
$n =50$ : $k =   49$ : $36x^2 + 18x - 1\, 801$	(Dress et Olivier, 1999, [3])\\
$n =100$ : $k =   90$ : $41x^2 + 33x - 43\, 321$	(Boston et Greenwood, 1995, \cite{Bos})\\
$n =200$ : $k = 166$ : $9x^2 + 3x - 16\, 229$		(Dress et Olivier)\\
$n =500$ : $k = 369$ : $x^2 + x - 1\, 354\, 363$		(Dress et Olivier)\\
$n = 1\,000$ : $k = 698$ : $x^2 + x - 1\, 354\, 363$		(Dress et Olivier)\\
\\[1ex]
Signalons \'egalement que les records $k$ sur $n$ pour le degr\'e 2 ont \'et\'e donn\'es
en continu jusqu'\`a $n = 40\, 338$ (Dress et Olivier, 1999, \cite{Dre})

	{\bf Degr\'e 3}
\\[1ex]
$n =50$ : $k=48$ : $x^3 + x^2 + 66x - 457$	(Dress et Olivier, 1999, \cite{Dre})\\
\hbox{}\hspace{2cm} ult\'erieurement 9 autres polyn\^omes avec 48 sur 50 ont \'et\'e trouv\'es \\
$n =100$ : $k=87$ : $3x^3 + x^2 - 17\, 888x - 365\, 413$\\
\hbox{}\hspace{2cm} (Dress et Landreau, pour ce record et tous ceux qui suivent)\\
$n =200$ : $k=154$ : $x^3 + x^2 + 2\, 764x + 16\, 553$ et $2x^3 + 2x^2 -10\, 664x -163\, 753$\\
$n =500$ : $k=338$ : $2x^3 + x^2 -13\, 145x - 218\, 651$\\
$n =1\, 000$ : $k = 601$ : $x^3 + x^2 - 48\, 610x + 4\, 021$

	{\bf Degr\'e 4}
\\[1ex]
$n =50$ : $k=49$ : 21 polyn\^omes\\
\hbox{}\hspace{2cm} (7 \`a coefficients entiers + 14 \`a coefficients rationnels non entiers)
			celui qui fournit le plus petit maximum des valeurs absolues prises est :
						$\frac14x^4 + \frac12 x^3 - \frac{537}4 x^2 + \frac{2\,459}2 x - 1\, 427$\\[1ex]
$n =100$ : $k =   88$ :  $\frac14x^4 + \frac12 x^3 - \frac{537}4 x^2 +\frac{1079}2 x  - 317$ \\[1ex]
$n =200$ : $k = 156$ :  $\frac34 x^4 + \frac12 x^3 - \frac{7511}4 x^2 +  \frac{7\,097}2x + 255\, 469$\\[1ex]
$n =500$ : $k = 316$ : $ \frac14 x^4+\frac12x^3-\frac{15\,921}4x^2+\frac{216\,359}2x+7\,422\,691$\\[1ex]
$n =1\,000$ : $k = 539$ :  $\frac14 x^4 + \frac12  x^3 + \frac{4303}4 x^2 + \frac{60\,395}2 x  - 2\, 092\, 807$\\

	{\bf Degr\'e 5}
\\[1ex]
$n =50$ : $k =   50$ (le record "$n$ sur $n$" est \`a 57)

	{\bf Degr\'e 6}
\\[1ex]
$n =50$ : $k =   50$ (le record "$n$ sur $n$" est \`a 58)

\subsection{M\'ethode utilis\'ee pour l'exploration num\'erique} 

	La base de la m\'ethode est une exploration exhaustive des
        polyn\^omes (normalis\'es) par boucles embo\^{\i}t\'ees sur
        les coefficients avec des limites fix\'ees, exploration
        convenablement cribl\'ee pour \'eviter les "petits" diviseurs
        des valeurs prises. La primalit\'e des valeurs est test\'ee
        par r\'ef\'erence \`a une table pr\'ealablement construite qui
        peut indexer jusqu'\`a 600 millions de nombres premiers.
 
	Il nous est apparu que la meilleure utilisation du crible sur
        les d.p.p. pour acc\'el\'erer la recherche consistait, dans
        son principe, \`a travailler, dans l'avant-derni\`ere boucle
        de l'exploration, sur les polyn\^omes sans terme constant
        $P^*(x) = a_0x_d + a_1x_{d-1} +\cdots + a_{d-1}x$.
 Pour un tel polyn\^ome, et pour un $p\le n$, 
nous d\'eterminons les classes o\`u doit se trouver $a_d$ modulo $p$
pour que les polyn\^omes $P(x) = P^*(x) + a_d$
 n'aient pas $p$ pour d.p.p. M\^eme en limitant ce premier crible \`a
 de "petits" $p$, le gain est tr\`es important.

	Pour le probl\`eme "$n$ sur $n$ ", ce crible fonctionne
        seul. Pour le probl\`eme "$k$ sur $n$", on lui adjoint un
        crible sur les d.p.p. "globalis\'e" comme nous l'avons
        expliqu\'e plus haut, en utilisant comme crit\`ere (empirique)
        une majoration de la somme des $\omega(p)/p$ pour $p$
        inf\'erieur \`a une borne convenable.

	Enfin, nous avons pris la d\'ecision de limiter la recherche au proche voisinage de l'origine, ce qui fait gagner un temps consid\'erable et ne laisse \'eventuellement \'echapper qu'une tr\`es faible proportion de "bons" polyn\^omes.

	La limite exacte des coefficients des polyn\^omes
        syst\'ematiquement explor\'es d\'epend du premier coefficient
        $a_0$, mais les zones d'exploration 
des coefficients ne sont
        pas homoth\'etiques, parce que le facteur de taille
        entra\^{\i}nerait alors une d\'ecroissance brutale des chances
        de succ\`es (cela est expliqu\'e plus loin dans la
        pr\'esentation du mod\`ele heuristique). Il serait fastidieux 
de d\'etailler le champ pr\'ecis des zones explor\'ees mais on peut donner des ordres de grandeur.
A chaque fois, le nombre de  polyn\^omes est consid\'er\'e avant de
faire op\'erer les diff\'erents cribles.
\begin{itemize}
\item degr\'e 2\\
\hbox{}\hspace{1cm}	$a_0$ de 1 \`a 3\,500, environ $10^{17}$ polyn\^omes consid\'er\'es,
\item degr\'e 3\\
\hbox{}\hspace{1cm} $a_0$ de 1 \`a 500, environ $10^{16}$ polyn\^omes consid\'er\'es, 
\item degr\'e 4  \\
\hbox{}\hspace{1cm}  - en coefficients entiers : $a_0$ de  1 \`a  36,
environ $3\,\cdot 10^{18}$ polyn\^omes,\\
\hbox{}\hspace{1cm} - en coefficients rationnels :
$b_0$ de 1 \`a  60, environ $1.4\,\cdot 10^{18}$ polyn\^omes,
\item degr\'e 5 \\ 
\hbox{}\hspace{1cm}  - en coefficients entiers :
$a_0$ de 1 \`a 5, environ $8 \,\cdot10^{17}$ polyn\^omes,\\
\hbox{}\hspace{1cm}  - en coefficients rationnels :
$b_0$ de 1 \`a 150, environ $5.2\,\cdot 10^{22}$ polyn\^omes,
\item degr\'e  6 \\ 
\hbox{}\hspace{1cm}  - en coefficients entiers :
$a_0$ de 1 \`a 2, environ $8 \,\cdot 10^{17}$ polyn\^omes,\\
\hbox{}\hspace{1cm}  - en coefficients rationnels :
$b_0$ de 1 \`a 20, environ $6.4\, \cdot 10^{23}$ polyn\^omes.

Nous avons fait tourner pour cela  des progammes \'ecrits en langage C et
parall\'elis\'es  principalement sur deux grappes de
calculs :\\
- la grappe Loire du CCIPL (Centre de Calcul Intensif des Pays de Loire,
http://www.cnrs-imn.fr/CCIPL/HTML/index.html),\\
- et la grappe Gaia du laboratoire LAGA  accessible depuis la plateforme de calcul du
Groupement De Services  Mathrice du CNRS (http://mathrice.org,  https://gaia.math.univ-paris13.fr).

\end{itemize}

\section{Le mod\`ele probabiliste heuristique dans le cas " $n$ sur
  $n$ "} 

\subsection{Localisation des zones riches en valeurs premi\`eres et
  normalisation} 

  	Comme on l'a d\'ej\`a not\'e, une heuristique \'el\'ementaire --- confirm\'ee par la recherche exp\'erimentale --- indique que le maximum de chances de trouver des valeurs premi\`eres pour un polyn\^ome $P$ consiste \`a les rechercher dans la zone des plus petites valeurs (absolues) prises. La traduction pour le mod\`ele de cette indication est compl\`etement diff\'erente en degr\'e 2 et en degr\'es sup\'erieurs.

	Dans le cas du degr\'e 2, la zone d'efficacit\'e peut se
        trouver tr\`es \'eloign\'ee de l'origine, c'est le cas
        notamment des polyn\^omes $a_0x^2 + a_1x + a_2$, avec $a_1$
        petit (normalisation de Boston et Greenwood) et $a_2$
        n\'egatif de grande valeur absolue. En cons\'equence, le
        mod\`ele probabiliste heuristique, pr\'esent\'e en \cite{Dre}
        pour le seul cas du degr\'e 2, "oubliait" enti\`erement les
        polyn\^omes eux-m\^emes et ne consid\'erait que des
        s\'equences de $n$ entiers \`a diff\'erences secondes (paires) constantes, caract\'eris\'ees par leurs valeurs extr\^emes et une valeur m\'ediane. 

En degr\'e sup\'erieur \`a 2 par contre, les essais que nous avons effectu\'es nous ont montr\'e que les polyn\^omes avec des petites valeurs loin de l'origine \'etaient statistiquement tr\`es rares. Comme il est beaucoup plus rapide de rechercher les valeurs premi\`eres au voisinage de l'origine que de chercher \`a localiser la zone o\`u $|F(x)|$ est "petit", nous avons donc suivi cette voie. De fa\c con plus pr\'ecise, nous avons constat\'e que le premier maximum relatif (maximum \'evoqu\'e au paragraphe 2.4 lors de la pr\'esentation du mur de Schinzel) du nombre de valeurs cons\'ecutives toutes premi\`eres se trouvait, \`a un nombre tr\`es faible d'exceptions pr\`es, dans la zone centrale tr\`es resserr\'ee entre $-n$ et $n$ (pour la variable). Cela a \'et\'e confirm\'e par les tests de recherche exhaustive que, pour $n\le 10$, nous avons pu effectuer jusqu'\`a la limite $n^d$. On peut \'egalement donner une justification heuristique de cette d\'ecroissance tr\`es rapide de la performance lorsqu'on s'\'ecarte de l'origine, et constater que la mod\'elisation est concordante avec les r\'esultats de l'exp\'erimentation.

	Il est bien s\^ur possible que la limitation au proche voisinage de l'origine fasse perdre quelques polyn\^omes "int\'eressants", mais le gain de rapidit\'e est tel qu'on regagne s\^urement beaucoup plus gr\^ace au nombre de polyn\^omes explor\'es. La situation est analogue dans le mod\`ele : il est trop compliqu\'e de calculer une estimation heuristique exacte du nombre des suites "int\'eressantes", mais les suites n\'eglig\'ees sont en proportion infime, et la minoration qui sera calcul\'ee sera extr\^emement voisine de la r\'ealit\'e.
	
	En d\'efinitive, le mod\`ele heuristique "$n$ sur $n$" en degr\'es sup\'erieurs \`a 2, pr\'esent\'e dans cet article, sera fond\'e sur la consid\'eration des polyn\^omes normalis\'es selon Boston et Greenwood, et en outre adapt\'e pour privil\'egier les "petites valeurs". Il est d\'ecrit ci-dessous.

\subsection{Pr\'esentation g\'en\'erale du mod\`ele} 

  	Le degr\'e $d$ \'etant fix\'e, $d\ge3$, on consid\`ere une
        famille croissante $(\F_B)$ de parties finies de l'espace
        $\Q^{d+1}$ des coefficients des polyn\^omes de degr\'e $d$, d\'ependant essentiellement d'un param\`etre de taille $B$  mais aussi de $n$.

	Les ensembles $\F_B$ seront \'equiprobabilis\'es. On
        \'etablira une estimation heuristique de la probabilit\'e
        $P_d(B,n)$ qu'un polyn\^ome pris au hasard dans $\F_B$ prenne
        une suite de $k$ valeurs cons\'ecutives premi\`eres, pour $k
        \ge n$ (i.e. $np^*(P)\ge  n$, avec la notation du paragraphe
        2.5).
 Cette probabilit\'e  $P_d(B,n)$ sera calcul\'ee \`a partir de la probabilit\'e $P_d(B,n,m)$  que les $n$ valeurs $P(m), P(m+1),\ldots, P(m+n-1)$ soient premi\`eres, probabilit\'e elle-m\^eme estim\'ee comme le produit d'un facteur "arithm\'etique" $C_d(n)$ par un facteur $T_d(B,n, m)$ qui traduit l'effet de taille du polyn\^ome $P$ consid\'er\'e.
	On peut alors \'evaluer la probabilit\'e $P_{d}(B,n)$, ainsi que le volume $N_d(B,n)$ de $\F_B$. L'indicateur crucial est l'esp\'erance math\'ematique heuristique $E_{d}(B,n) =
P_{d}(B,n)N_d(B,n)$ du nombre de polyn\^omes $P$ de $\F_B$ qui satisfont $np^*(P)\ge  n$, et on termine l'\'etude du mod\`ele en examinant  le comportement de cette esp\'erance lorsque $B$ tend vers l'infini.

Les parties $\F_B$ seront des pav\'es de l'espace des coefficients,
qui seront param\'etris\'es de telle sorte que $B$ majore le maximum de la valeur absolue des valeurs prises par les polyn\^omes de $\F_B$
sur l'intervalle $[-n, n]$. 
Bien entendu, il n'y pas a priori de mani\`ere unique de d\'efinir les $\F_B$, on donnera au paragraphe 4.5 ceux que nous avons utilis\'es, pour lesquels les diff\'erentes bornes impos\'ees aux coefficients sont d\'etermin\'ees \`a partir de consid\'erations empiriques pour s'adapter aux (relativement) fortes probabilit\'es de primalit\'e au tr\`es proche voisinage de l'origine.

\subsection{S\'eparation du facteur arithm\'etique et du facteur de taille} 

 	Etant donn\'e un polyn\^ome g\'en\'erique, i.e. tir\'e au hasard avec \'equiprobabililit\'e sur le pav\'e $\F_B$ de l'espace des coefficients, il faut commencer par donner une estimation heuristique de la probabilit\'e $P_{d}(B,n, m)$ que les $n$ valeurs $P(m), P(m+1),\ldots, P(m+n-1)$ soient premi\`eres. Cela se fait par un criblage par les nombres premiers $p$ inf\'erieurs \`a la racine carr\'ee des valeurs (absolues) prises, que l'on applique de fa\c con diff\'erente pour les valeurs de $p$ inf\'erieures \`a $n$ et pour celles sup\'erieures.
	Pour un $p\le n$, la p\'eriodicit\'e rend non ind\'ependantes les divisibilit\'es par $p$ des valeurs de $P(m)$ \`a $P(m+n-1$), et la condition de non-divisibilit\'e est globale : que $p$ ne soit pas d.p.p. du polyn\^ome (g\'en\'erique) $P$. Pour $n$ fix\'e et $B$ grand, cette probabilit\'e ne d\'epend pas de $m$, et l'on obtient ainsi pour chaque $p$ un facteur $D_d(p)$. On suppose ensuite que, dans le mod\`ele, les divisibilit\'es par les nombres premiers successifs sont ind\'ependantes, ce qui est une approximation raisonnable (en effet dans le cas de tous les nombres entiers, il n'y a pas plus qu'un facteur $2e^{\gamma}= 1.123\ldots$ entre le th\'eor\`eme de Mertens et l'heuristique des divisibilit\'es ind\'ependantes). Cela permet alors de multiplier les $D_d(p)$ pour $p\le n$, obtenant ainsi le facteur "arithm\'etique" (g\'en\'erique) $C_d(n)$ annonc\'e.
	Par contre, pour les $p > n$, on crible chaque valeur de $P(m)$ \`a $P(m+n-1)$ individuellement. Les valeurs de $m$ interviennent "\`a travers" la taille du polyn\^ome, dont le facteur obtenu $T_d(B,n, m)$ traduit l'effet (de fa\c con assez complexe comme on le verra).

	Cette mod\'elisation conduit \`a la formule
$$
P_d(B,n, m) = \prod_{p\le n} D_d(p)\prod_{x=m}^{m+n-1}\prod_{n< p\le \sqrt{|P(x)|}}(1-\frac1p).
$$ 
Les facteurs $\prod_{n< p\le\sqrt{|P(x)|}}(1-\frac1p)$ sont
estim\'es avec le th\'eor\`eme de Mertens par $\frac{\log n}{\log\sqrt{|P(x)|}}$. On a alors la d\'ecomposition
$$
P_d(B,n, m) \sim C_d(n) T_d(B,n, m),
$$
avec $C_d(n) =\prod_{p\le n} D_d(p)$  et  $T_d(B,n, m) = \prod_{x=m}^{m+n-1}\frac{\log n}{\log\sqrt{|P(x)|}} $.

\subsection{Le facteur arithm\'etique}

	Il faut commencer par estimer chaque facteur $D_d(p)$. Nous le
        ferons en nous limitant au cas des coefficients entiers. Comme
        $n$ sera fix\'e et que $B$ tendra vers l'infini, le tirage au
        hasard du polyn\^ome g\'en\'erique de degr\'e $d$, avec \'equiprobabililit\'e
        sur le pav\'e $\F_B$ de l'espace des coefficients, sera
        convenablement mod\'elis\'e par l'\'equiprobabilit\'e des
        coefficients modulo $p$ sur le produit  ${\mathbb F}_p^{d+1}$, (i.e. une probabilit\'e uniforme \'egale \`a $\frac1{p^{d+1}}$  pour tout $(d+1)$-uplet de coefficients modulo $p$). On n\'eglige les polyn\^omes dont une valeur serait exactement $p$ ou $-p$, ce qui a une incidence quantitative infime. Cette mod\'elisation, justifi\'ee par l'\'equivalence asymptotique, conservera une justification heuristique forte alors m\^eme qu'elle sera utilis\'ee comme nous le ferons pour des valeurs fix\'ees de $n$ et de $B$.

\begin{prop}
Dans le cas des polyn\^omes \`a coefficients entiers de degr\'e  $d$, la probabilit\'e
heuristique que $p$ ne soit pas d.p.p. vaut :
 $$D_d(p) = (1-\frac1p)\sum_{k=0}^{\min(d,p-1)} (-1)^k \binom{p-1}{k} p^{-k}.$$
En particulier, pour $p\le d+1$ , on a $D_d(p) = (1-\frac1p)^p$.
\end{prop}
\noindent D\'emonstration: il s'agit de d\'eterminer le cardinal $M_d(p)$ de l'ensemble des polyn\^omes d\'efinis modulo $p$ de degr\'e inf\'erieur ou \'egal \`a $d$  et sans z\'ero dans $\FF_p$.
\\[1ex]
On travaille n\'ecessairement dans l'ensemble des polyn\^omes de terme constant non nul. Consid\'erons  parmi ceux-ci uniquemement les polyn\^omes de terme constant   \'egal \`a 1, il suffira ensuite  de multiplier le cardinal obtenu  par $(p-1)$. 
	Pour tout $i$ de 1 \`a $p-1$, on pose $E_i = \{F \in {\FF}_p[X]\; |\;
         \deg F \le d, F(0)=1, F(i) = 0\}$. Un calcul rapide montre
         que $|E_i| = p^{d-1}$ et on peut de m\^eme calculer, pour $i_1 < i_2
          <\ldots < i_k$, $1 \le k \le \min(d,p-1)$, $|E_{i_1} \cap E_{i_2}\cap \ldots \cap E_{i_k}| = p^{d-k}$.
 Le  principe combinatoire classique d'inclusion-exclusion montre alors
 que le nombre de  polyn\^omes sans z\'ero dans $\FF_p$ et de terme constant 1 vaut
$$\sum_{k=0}^{m} (-1)^k \binom{p-1}{k} p^{d-k},$$
o\`u $m=\min(d,p-1)$ car les polyn\^omes concern\'es ne peuvent avoir plus de $\min(d,p-1)$ z\'eros sur $\FF_p$.
On en d\'eduit la valeur de $M_d(p)$ par multiplication par $(p-1)$,
puis la valeur  de $D_d(p)$ :
$$D_d(p) =M_d(p)/p^{d+1}=(1-\frac1p)\sum_{k=0}^{m} (-1)^k \binom{p-1}{k} p^{-k}.$$
On remarque que pour $p\le d+1$, on obtient la formule simplifi\'ee 
$D_d(p)=(1-\frac1p)^p$.
\begin{est}
 On donne dans le tableau suivant les valeurs  num\'eriques du facteur arithm\'etique
$C_d(n)\cdot 10^{6} $, pour $d= 3, 4, 5, 6$ et $n = 20, 25, \ldots, 50$
ainsi que la valeur de $\prod_{p\le n}(1-\frac1p)^p \cdot 10^{6}$.
 $$
 \begin{array}{|c|c|c|c|c|c|c|c|}
\hline
 d \backslash n & 20 & 25 & 30 & 35 & 40 & 45 &50 \\
\hline
3 &106.7 
& 36.10 
& 12.19
& 4.111
& 1.385 
& 0.1569 
& 0.0528 
\\
\hline 
4 & 133.8 
& 48.60 
& 17.76 
& 6.503 
&  2.390
& 0.3242 
& 0.1196 
\\
\hline 
5 & 130.1 
& 46.76 
&  16.88
& 6.097 
& 2.208 
& 0.2906 
&  0.1055
\\
\hline 
6 & 130.5 
& 46.95 
& 16.97 
&  6.143
& 2.230 
&  0.2947
&  0.1073 
\\
\hline
\prod_{p\le n}(1-1/p)^p & 130.5
& 46.94
& 16.96
& 6.139
& 2.227
& 0.2943
& 0.1071
\\
\hline
\end{array}
$$
\end{est}
\tabcaption{valeurs  num\'eriques du facteur arithm\'etique
  $C_d(n)\cdot 10^{6}$}

	Ce tableau est obtenu par calcul num\'erique de $C_d(n) =\prod_{p\le
          n}  D_d(p)$ avec les valeurs de $D_d(p)$ de la proposition
        1.
On remarquera qu'\`a partir de $d\ge5$ les valeurs de $C_d(n)$ sont
tr\`es proches de la valeur $\prod_{p\le n}(1-1/p)^p$ ind\'ependante
de $d$. La raison en est que  la contribution principale est en fait celle des nombres premiers $p\le d+1$ pour laquelle on dispose de  la formule simplifi\'ee.
D'un point de vue asymptotique, on peut facilement montrer que 
$$ \prod_{p\le n}(1-1/p)^p\sim Ke^{-\pi(n)}/\sqrt{\log(n)} \quad \mbox{
  o\`u }  K\approx 0.7.
$$
      
\subsection{Le facteur de taille}

	Le facteur de taille $T_d(B,n,m) =  \prod_{x=m}^{m+n-1}\frac{\log n}{\log\sqrt{|P(x)|}}$ a \'et\'e introduit au
        d\'ebut du paragraphe 4.3. Nous avions annonc\'e au paragraphe
        4.2 une d\'efinition des pav\'es $\F_B$  de l'espace des
        coefficients qui ferait de $B$ un majorant du maximum de la
        valeur absolue des valeurs prises par les polyn\^omes de
        $\F_B$ sur l'intervalle $[-n, n]$. N\'eanmoins, dans
        l'estimation du facteur de taille, ce serait une maladresse de
        consid\'erer seulement cette borne $B$ et de minorer
        $T_d(B,n,m)$ par $\left(\frac{\log n}{\log\sqrt{B}}\right)^n$. 
En proc\'edant ainsi, non seulement on obtiendrait une minoration trop inf\'erieure \`a une estimation raisonnable (et d'ailleurs en d\'esaccord avec les r\'esultats num\'eriques), mais surtout on cacherait compl\`etement le resserrement de la zone d'efficacit\'e autour de l'origine. Il faut donc choisir la "forme" efficace des parties $\F_B$ et la faire intervenir  dans les estimations.

	Nous avons 
consid\'er\'e des parties $\F_B$ 
ensembles des $(d+1)$-uplets des coefficients des polyn\^omes dont les coefficients $|a_j|$ sont major\'es par des bornes plus petites que $ \frac{B}{(d+1)n^{d-j}}$, et ce d'autant plus que le degr\'e $d-j$ est faible. Ce crit\`ere qualitatif garantit (heuristiquement et exp\'erimentalement) que les parties $\F_B$ contiennent la zone d'efficacit\'e o\`u se trouvent la plupart des "bons" polyn\^omes, et les r\'esultats ne d\'ependent que faiblement du choix pr\'ecis --- math\'ematiquement arbitraire. Notre choix s'appuie sur la limitation $|a_d| \le B^{3/4}$ pour le dernier coefficient, limitation tr\`es voisine de celle adopt\'ee dans nos "explorations" exp\'erimentales.  Nous avons ainsi pris :
\begin{itemize}
\item  	$1 \le a_0 \le \frac{B}{n^d}$, 
\item	$0 \le a_1 \le \frac{d}2 a_0$  (impos\'e par la normalisation de Bonston et Greenwood),
\item $	|a_j| \le \frac{B^{1-\frac{j}{4d}}}{n^{d-j}}$ pour $2\le j\le d$.
\end{itemize}

(Pour simplifier les formules, nous continuerons \`a nous limiter au
cas des polyn\^omes \`a coefficients entiers.) Il faut imm\'ediatement
noter la contrainte $B\ge n^d$ (n\'ecessaire pour que $a_0$ puisse au moins prendre la valeur 1), apparemment triviale, mais qui en fait jouera un r\^ole important quand on voudra effectuer des pr\'edictions.

	Pour calculer la probabilit\'e ponctuelle, nous \'ecrirons 
$T_d(B,n, m) = \left(\frac{\log n}{\log\sqrt{B}}\right)^n K_d(B,n, m)$, et
nous \'etablirons une estimation efficace du facteur correctif $K_d(B,n,
m)=\prod_{x=m}^{m+n-1}\frac{\log B}{\log |P(x)|}$, pris pour un polyn\^ome "typique" du pav\'e $\F_B$.

	Il faudra ensuite passer de l'estimation de la probabilit\'e
        ponctuelle $P_d(B,n,m)$ \`a la probabilit\'e globale
        $P_d(B,n)$ qu'un polyn\^ome g\'en\'erique pr\'esente (au
        moins) $n$ valeurs cons\'ecutives premi\`eres. Comme il s'agit
        d'\'ev\`enements (tr\`es) rares, on peut assimiler les
        probabilit\'es aux esp\'erances, ce qui permet d'approcher la
        probabilit\'e de la r\'eunion d'\'ev\`enements non
        ind\'ependants par la somme de leurs probabilit\'es. On
        obtient ainsi  $P_d(B,n)=\sum_{m\in\Z} P_d(B,n, m)$, et l'estimation de cette somme s'effectuera de fa\c con standard par assimilation \`a une int\'egrale.

{\bf a}) Estimation du facteur \'el\'ementaire du facteur correctif $K_d(B,n, m)$.\\[1ex]
Pour simplifier la suite de la pr\'esentation et des calculs, on pose
$n' = \left[\frac{n}2\right]$; compte tenu des nombreuses
approximations effectu\'ees, on n'introduit aucune erreur
significative en effectuant les calculs comme si on avait exactement
$n'=\frac{n}2$ et comme si l'intervalle $[-n', n']$ contenait
exactement $n$ entiers.

	C'est maintenant qu'intervient la d\'efinition pr\'ecise
        des parties $\F_B$ qui vient d'\^etre effectu\'ee : elle
        implique que le polyn\^ome g\'en\'erique $P$ prend \`a
        l'origine des valeurs de l'ordre de $B^{3/4}$, que son terme
        de plus haut degr\'e est de l'ordre de  $\frac{B}{n^d}x^d$, et
        que ce terme domine num\'eriquement tous les autres sauf dans
        un voisinage restreint de l'origine. Pour donner une
        estimation analytique du facteur correctif, nous effectuerons
        la sommation des valeurs de $l(x) := \log\left( \frac{\log
            B}{\log |P(x)|}\right)$, cette fonction \'etant
        elle-m\^eme \'evalu\'ee par son approximation sym\'etrique non
        triviale la plus simple, i.e. $h(x) = a + bx^2$. Cela
        r\'ealise un compromis efficace entre une approximation
        simpliste et des calculs impraticables.

	Compte tenu des indications donn\'ees ci-dessus sur l'ordre de
        grandeur du  polyn\^ome g\'en\'erique $P$, cette approximation
        sera ajust\'ee pour $x=0$ et $x=\pm n'$ :
        \begin{itemize}
        \item[-] en $x=0$ : $h(0)\approx l(0)= \log\left( \frac{\log
            B}{\log B^{3/4}}\right)=\log \frac43$,
\item[-] en $x=n'$ : $h(n')\approx l(n')=  \log\left( \frac{\log B}{\log
      \frac{B}{2^d}}\right)
\approx \frac{d\log 2}{\log B}$.
\end{itemize}
On trouve alors $a =\log \frac43$ et $b=-\frac4{n^2}c \log \frac43$ en ayant
pos\'e 
$c=1-\frac{d\log 2}{\log B\log \frac43}$.

	Cette approximation n'est raisonnable que si la valeur
        $|P(n')|=\frac{B}{2^d}$  prise en assimilant le polyn\^ome \`a
        son terme de plus haut degr\'e est sup\'erieure \`a l'ordre de
        grandeur $B^{3/4}$ du terme constant. On devra donc respecter
        la contrainte $B\ge 2^{4d}$, ce qui d\'ecoule, si $n\ge 16$,
        de la contrainte plus forte $B\ge n^d$ \'evoqu\'ee juste
        apr\`es la d\'efinition des $\F_B$. On peut enfin noter que
        l'in\'egalit\'e $c>0$, n\'ecessaire pour que l'origine soit un
        maximum de $h(x)$, est alors trivialement satisfaite.

	On peut enfin comparer les valeurs $l(n) = 0$ et 
$h(n) = \frac{4d\log 2 }{\log B}-3\log\frac43$ qui est n\'egatif si
$B\ge n^d$ et $n\ge25$ : notre approximation \`a partir de la fonction $h$ conduira donc \`a une sous-estimation de la probabilit\'e $P_d(B,n,m)$ lorsque l'on s'\'ecartera de l'intervalle central, et par cons\'equent \`a une sous-estimation globale, mais ce n'est pas vraiment g\^enant car $P_d(B,n, m)$ d\'ecroit rapidement et les "queues de distribution" sont sans incidence num\'erique notable.

{\bf b)} Fin de l'estimation du facteur correctif $K_d(B,n,m)$
\\[1ex]
On approche classiquement $\log K_d(B,n,m)\approx \sum_{x=m}^{m+n-1} h(x)$
  par l'int\'egrale de $h$ correspondante, soit
  \begin{equation*}
    \begin{split}
\log K_d(B,n,m)&\approx \int_{m}^{m+n} h(x) dx = an +\frac{b}3\left((m+n)^3-m^3\right)\\
&=\log\frac43\left(n-\frac4{n^2}c(m^2n+mn^2+n^3/3)\right)\\
 &= - \log\frac43\left(\frac{4c}n(m^2+mn)+n(\frac{4c}3-1)\right)
=- \log\frac43\left(\frac{4c}n(m+n/2)^2+n(\frac{c}3-1)\right).
\end{split}
\end{equation*}

{\bf c}) Estimation de la probabilit\'e globale $P_d(B,n)$.
\\[1ex]
On termine en calculant comme annonc\'e plus haut la somme 
$P_d(B,n) =\sum_{m\in\Z} P_d(B,n, m)$, o\`u $P_d(B,n, m)\approx
C_d(n) T_d(B,n, m)$, avec 
$T_d(B,n,m) = \left(\frac{\log n}{\log \sqrt{B}}\right)^n K_d(B,n, m)$.

On est donc ramen\'e \`a \'evaluer $\sum_{m\in\Z} K_d(B,n,m)$,  ce qui se fait en l'approchant classiquement par une int\'egrale :
 \begin{equation*}
   \begin{split}
   \sum_{m\in\Z}  K_d(B,n,m) &\approx \exp(-n \log\frac43(\frac{c}3-1)) 
\int_{-\infty}^{+\infty} \exp(-4c\frac{\log\frac43}{n}(x+\frac{n}2)^2) dx\\
&= \left(\frac43\right)^{n(1-\frac{c}3)}
\sqrt{\frac{n}{4c\log\frac43}}\int_{-\infty}^{+\infty}
 \exp(-u^2) du\\
&=\left(\frac43\right)^{n(1-\frac{c}3)}\sqrt{\frac{\pi n}{4c\log\frac43}}.
\end{split}
\end{equation*}

Lorsque $n$ et $d$ sont fix\'es et que $B$ tend vers
l'infini, $c=1-\frac{d\log 2}{\log B\log \frac43}$ tend (lentement) vers 1, et par
cons\'equent $\sum_{m\in\Z} K_d(B,n,m)$  est \'equivalent \`a 
$ \left(\frac43\right)^{2n/3}\sqrt{\frac{\pi n}{4\log\frac43}}$ . Mais on conservera $c$ dans l'expression de la proposition 2 ci-dessous, de fa\c con \`a pouvoir la remplacer par sa valeur exacte dans les tests de comparaison entre le mod\`ele heuristisque et les exp\'erimentations num\'eriques qui seront effectu\'es au paragraphe suivant, et \`a conserver la pertinence de ces comparaisons.   

\begin{prop} La probabilit\'e globale $P_d(B,n)$ peut \^etre
  estim\'ee de la fa\c con suivante lorsque, $d$ et n \'etant fix\'es, $B$ tend vers l'infini 
$$
P_d(B,n) \approx C_d(n) \left(\frac{\log n}{\log\sqrt{B}}\right)^n
 \left(\frac43\right)^{n(1-\frac{c}3)}\sqrt{\frac{\pi n}{4c\log\frac43}}.
$$
\end{prop}

D\'emonstration: on a
\begin{equation*}
  \begin{split}
		P_d(B,n) &\approx \sum_{m\in\Z} P_d(B,n, m)
                \approx C_d(n) \sum_{m\in\Z} T_d(B,n, m)\\
&\approx C_d(n)\left(\frac{\log n}{\log\sqrt{B}}\right)^n\sum_{m\in\Z} K_d(B,n, m)\\
&\approx C_d(n)\left(\frac{\log n}{\log\sqrt{B}}\right)^n
\left(\frac43\right)^{n(1-\frac{c}3)}
\sqrt{\frac{\pi n}{4c\log\frac43}}.
\end{split}
\end{equation*}

\subsection{Etude de l'esp\'erance globale, mod\`ele heuristique} 

	Nous pouvons maintenant conclure. Il reste \`a calculer le nombre de points $N_d(B,n)$ d'un ensemble $\F_B$ et le produit de la probabilit\'e globale $P_d(B,n)$ par $N_d(B,n)$ fournit une estimation de l'esp\'erance math\'ematique $E_d(B,n)$ du nombre de polyn\^omes pour lesquels $np^*(P)\ge n$.

\begin{prop} Le nombre de points $N_d(B,n)$ des parties $\F_B$
  d\'efinies ci-dessus peut \^etre estim\'e par :
$$
N_d(B,n) \approx  d
2^{d-3}\frac{B^{\frac{7d}8+\frac78+\frac1{4d}}}{n^{\frac{d^2}2+\frac{d}2+1}}.
$$
\end{prop}

D\'emonstration: Compte-tenu des contraintes $1 \le a_0\le
\frac{B}{n^d}$  et $0 \le a1 \le \frac{d}2a_0$, le nombre de valeurs
possibles des couples $(a_0,a_1)$ peut \^etre estim\'e par
$\frac{dB^2}{4n^{2d}}$; par contre les contraintes sur les
coefficients $a_j, 2 \le j \le n,$ sont ind\'ependantes de $a_0$, et
$N_d(B,n)$ peut donc \^etre estim\'e par le produit
$$
\frac{dB^2}{4n^{2d}}\times \frac{2B^{1-\frac2{4d}}}{n^{d-2}}\times\ldots
\times \frac{2B^{1-\frac{d}{4d}}}{n^{d-d}}.
$$
	Nous pouvons maintenant conclure. Le produit de la probabilit\'e globale $P_d(B,n)$ par le nombre de points $N_d(B,n)$ fournit une estimation de l'esp\'erance math\'ematique $E_d(B,n)$ du nombre de polyn\^omes pour lesquels $np^*(P)\le n$.
      
      \begin{prop} L'esp\'erance math\'ematique $E_d(B,n)$ d\'efinie
          ci-dessus peut \^etre estim\'ee de la fa\c con suivante
          lorsque, $d$ et $n$ \'etant fix\'es, $B$ tend vers l'infini :
$$
E_d(B,n) \approx H_d(n) c^{-1/2} \left(\frac43\right)^{n(1-\frac{c}3)}
(\log B)^{-n} B^{\frac{7d}8+\frac78+\frac1{4d}},
$$
avec
$H_d(n) = \sqrt{\frac{\pi}{\log\frac43}} \,d\, C_d(n) 2^{n+d-4} (\log n)^n
n^{-\frac{d^2}2-\frac{d}2-\frac12}$.
\end{prop}

D\'emonstration: On effectue comme annonc\'e le produit de $P_d(B,n)$
(proposition 2) par $N_d(B,n)$ (proposition 3), et on isole le facteur
$c^{-1/2} \left(\frac43\right)^{n(1-\frac{c}3)}
B^{\frac{7d}8+\frac78+\frac1{4d}} (\log B)^{-n}$   
 qui regroupe les facteurs \'el\'ementaires qui d\'ependent de $B$
(on rappelle que $c=1-\frac{d\log 2}{\log B\log \frac43}$). Il reste
alors le facteur $H_d(n)$ qui regroupe les facteurs \'el\'ementaires
qui ne d\'ependent que de $d$ et de n, comme d\'etaill\'e dans la proposition 5.

	Une analyse fine de la grandeur et de la rapidit\'e de croissance/d\'ecroissance des diff\'erents facteurs de $E_d(B,n)$ sera effectu\'ee \`a la fin du paragraphe 4.7. (ils ont \'et\'e \'ecrits dans l'\'enonc\'e de la proposition 4 par ordre croissant de vitesse de variation). Mais on peut noter d\`es maintenant que la croissance du facteur $B^{7d/8}$ est fortement pr\'epond\'erante et que par cons\'equent (et comme on pouvait s'y attendre), $d$ et $n$ \'etant fix\'es, $E_d(B,n)$ tend vers l'infini avec $B$.

	Cette formule asymptotique n'a qu'une valeur indicative. Pour
        comparer les valeurs num\'eriques heuristiques avec les
        r\'esultats exp\'erimentaux, nous prendrons les formules
        "compliqu\'ees" de mani\`ere \`a ne pas introduire de
        "brouillage" num\'erique inutile.
      
\subsection{Etude de l'esp\'erance globale et validation exp\'erimentale}

	Dans la comparaison des valeurs num\'eriques du mod\`ele heuristique avec les r\'esultats de nos exp\'erimentations, on pr\'esente, en plus des esp\'erances, les nombres de polyn\^omes dans $\F_B$ et les probabilit\'es, de fa\c con \`a montrer clairement leurs ordres de grandeur respectifs. Les valeurs de $N_d(B,n)$ sont les valeurs exactes (les valeurs de la proposition 3 souffrant de remplacement de parties enti\`eres par des valeurs r\'eelles, ce qui est tr\`es mauvais lorsque $B$ est de  l'ordre de $n^d$) et  les probabilit\'es sont calcul\'ees 
avec les vraies valeurs de $C_d(n)$  et
  la vraie valeur de $c$. Enfin, le nombre r\'eel $R_d(B,n)$ de
polyn\^omes dans le pav\'e
$\F_B$  avec au moins  $n$ sur $n$  valeurs premi\`eres est le nombre trouv\'e exp\'erimentalement sans la contrainte que les valeurs premi\`eres soient distinctes en valeur absolue (contrainte convenue par Boston et Greenwood pour les records, mais artificielle et surtout impossible \`a respecter dans les calculs heuristiques de probabilit\'es).

	On rappelle que les polyn\^omes $P$ des pav\'es
$\F_B$  v\'erifient entre autres deux propri\'et\'es : sur l'intervalle $[-n,n]$ on a $|P(x)| \le B$, et le terme constant $a_d$ est inf\'erieur \`a $B^{3/4}$ en valeur absolue.
La borne param\'etrique $B$ est pr\'esent\'ee sous la forme $B= k n^d$ 
		et c'est la valeur de $k$  qui est port\'ee dans le
                tableau.
La derni\`ere ligne du tableau donne le quotient de $R =R_d(B, n)$ par $E = E_d(B, n)$, qui mesure l'ad\'equation du mod\`ele.

\noindent {\bf d = 3} :

\renewcommand{\arraystretch}{1.6} 
$n = 20$  ($C_3(20) = 1.067 \cdot 10^{-4}$)
$$ \scriptsize \begin{array}{|c|c|c|c|c|c|c|c|c|c|} \hline k     
& 1  & 2  & 4  & 8  & 16  & 32  & 64  & 128  & 256   
 \\ \hline     
  N_3(B,n)    
& 6.1\cdot 10^{5}  & 5.4\cdot 10^{6}  & 4.9\cdot 10^{7}  & 4.9\cdot 10^{8}  & 5.3\cdot 10^{9}  & 6\cdot 10^{10}  & 6.9\cdot 10^{11}  & 8.2\cdot 10^{12}  & 9.7\cdot 10^{13}  
 \\ \hline   
   P_3(B,n)   
& 1.2\cdot 10^{-4}  & 2.1\cdot 10^{-5}  & 4.3\cdot 10^{-6}  & 1\cdot 10^{-6}  & 2.7\cdot 10^{-7}  & 7.6\cdot 10^{-8}  & 2.3\cdot 10^{-8}  & 7.7\cdot 10^{-9}  & 2.7\cdot 10^{-9}  
  \\ \hline     
  E_3(B,n)   
 & 70.3  & 112.7  & 210.5  & 496.1  & 1400.8  & 4527.4  & 16230.4  & 63120.7  & 262563  
  \\ \hline   
  R_3(B,n)   
 & 97  & 297  & 875  & 2552  & 7516  & 23824  & 78362  & 267978  & 968007  
  \\ \hline    
   R/E   
 & 1.38  & 2.64  & 4.16  & 5.14  & 5.37  & 5.26  & 4.83  & 4.25  & 3.69  
 \\ \hline \end{array} $$
\tabcaption{ Comparaisons du  mod\`ele heuristique avec les valeurs
  exp\'erimentales pour  $d=3$, $n=20$}

\bul $n=30$, ($C_3(30)=1.219\cdot 10^{-5}$),

$$ \scriptsize \begin{array}{|c|c|c|c|c|c|c|c|c|c|} \hline k     
& 1  & 2  & 4  & 8  & 16  & 32  & 64  & 128  & 256   
 \\ \hline  
  N_3(B,n)    
& 2.8\cdot 10^{6}  & 2.5\cdot 10^{7}  & 2.2\cdot 10^{8}  & 2.2\cdot 10^{9}  & 2.4\cdot 10^{10}  & 2.7\cdot 10^{11}  & 3.2\cdot 10^{12}  & 3.7\cdot 10^{13}  & 4.5\cdot 10^{14}  
  \\ \hline    
   P_3(B,n)   
 & 2.6\cdot 10^{-6}  & 2.9\cdot 10^{-7}  & 3.9\cdot 10^{-8}  & 5.9\cdot 10^{-9}  & 10\cdot 10^{-10}  & 1.9\cdot 10^{-10}  & 3.8\cdot 10^{-11}  & 8.6\cdot 10^{-12}  & 2\cdot 10^{-12}  
  \\ \hline    
  E_3(B,n)  
& 7.1  & 7.3  & 8.7  & 13.2  & 24.1  & 51.2  & 122.1  & 320.6  & 911.6  
  \\ \hline    
   R_3(B,n)   
 & 13  & 32  & 61  & 111  & 178  & 329  & 621  & 1168  & 2412  
  \\ \hline    
  R/E    
& 1.82  & 4.39  & 7.01  & 8.43  & 7.39  & 6.43  & 5.08  & 3.64  & 2.65  
  \\ \hline \end{array} 
$$    
 \tabcaption{ Comparaisons du  mod\`ele heuristique avec les valeurs
  exp\'erimentales pour  $d=3$, $n=30$}

\bul $n=40$, ($C_3(40)= 1.385\cdot 10^{-6}$),
$$ \scriptsize \begin{array}{|c|c|c|c|c|c|c|c|c|c|} \hline k     
    & 4  & 8  & 16  & 32  & 64  & 128  & 256  & 512  & 1024  
  \\ \hline    
   N_3(B,n)   
 &6.6\cdot10^{8}  &6.6\cdot10^{9}  &7.1\cdot10^{10}  &8\cdot10^{11}  &9.3\cdot10^{12}  &1.1\cdot10^{14}  &1.3\cdot10^{15}  &1.6\cdot10^{16}  &1.9\cdot10^{17}  
  \\ \hline    
   P_3(B,n)   
 &3.6\cdot10^{-10}  &3.5\cdot10^{-11}  &3.9\cdot10^{-12}  &4.9\cdot10^{-13}  &6.9\cdot10^{-14}  &1\cdot10^{-14}  &1.7\cdot10^{-15}  &3.1\cdot10^{-16}  &6\cdot10^{-17}  
  \\ \hline    
   E_3(B,n)   
 & 0.2  & 0.2  & 0.3  & 0.4  & 0.6  & 1.2  & 2.3  & 4.9  & 11.3  
  \\ \hline    
   R_3(B,n)   
 & 1  & 6  & 7  & 8  & 9  & 12  & 15  & 15  & 15  
  \\ \hline    
   R/E   
 & 4.24  & 25.9  & 25.1  & 20.2  & 14  & 10.4  & 6.56  & 3.06  & 1.33  
  \\ \hline \end{array} $$    
 \tabcaption{ Comparaisons du  mod\`ele heuristique avec les valeurs
  exp\'erimentales pour  $d=3$, $n=40$}

\noindent {\bf d = 4} :	

\bul $n=20$, ($C_4(20)=1.338\cdot 10^{-4}  $),
$$  \begin{array}{|c|c|c|c|c|c|} \hline k 
& 1  & 2  & 4  & 8  & 16   
 \\ \hline    
   N_4(B,n)   
 &1.5\cdot10^{10}  &2.1\cdot10^{11}  &3.4\cdot10^{12}  &6.2\cdot10^{13}  &1.2\cdot10^{15}  
  \\ \hline    
   P_4(B,n)   
 &4.6\cdot10^{-7}  &1.2\cdot10^{-7}  &3.7\cdot10^{-8}  &1.2\cdot10^{-8}  &4\cdot10^{-9}  
  \\ \hline    
   E_4(B,n)   
 & 6709.3  & 26206.8  & 125703  & 727070  & 4.9D+06  
  \\ \hline    
   R_4(B,n)   
 & 15838  & 83645  & 470634  & 2855234  & 18639074  
  \\ \hline    
   R/E   
 & 2.36  & 3.19  & 3.74  & 3.93  & 3.84  
  \\ \hline \end{array} 
$$
\tabcaption{ Comparaisons du  mod\`ele heuristique avec les valeurs
  exp\'erimentales pour  $d=4$, $n=20$}

\bul $n=30$, ($C_4(30)= 1.776\cdot 10^{-5}$),
 $$  \begin{array}{|c|c|c|c|c|c|} \hline k 
& 1  & 2  & 4  & 8  & 16   
 \\ \hline     
  N_4(B,n)    
&2.2\cdot10^{11}  &3.2\cdot10^{12}  &5.3\cdot10^{13}  &9.5\cdot10^{14}  &1.9\cdot10^{16}   
 \\ \hline     
  P_4(B,n)    
&6.7\cdot10^{-10}  &1.3\cdot10^{-10}  &2.7\cdot10^{-11}  &6.2\cdot10^{-12}  &1.5\cdot10^{-12}   
 \\ \hline     
  E_4(B,n)    
& 150.7  & 419.7  & 1438.4  & 5937  & 28425.9  
 \\ \hline     
  R_4(B,n)    
& 494  & 1700  & 6019  & 22062  & 89354   
 \\ \hline   
  R/E    
& 3.28  & 4.05  & 4.18  & 3.72  & 3.14  
 \\ \hline \end{array} $$    
 \tabcaption{ Comparaisons du  mod\`ele heuristique avec les valeurs
  exp\'erimentales pour  $d=4$, $n=30$}

\bul $n=40$, ($C_4(40)=2.390 \cdot 10^{-6}$),
$$  \begin{array}{|c|c|c|c|c|c|} \hline k  
& 1  & 2  & 4  & 8  & 16  
 \\ \hline     
  N_4(B,n)    
&1.6\cdot10^{12}  &2.3\cdot10^{13}  &3.7\cdot10^{14}  &6.6\cdot10^{15}  &1.3\cdot10^{17}   
 \\ \hline     
  P_4(B,n)    
&1\cdot10^{-12}  &1.4\cdot10^{-13}  &2.1\cdot10^{-14}  &3.4\cdot10^{-15}  &6.1\cdot10^{-16}   
 \\ \hline     
  E_4(B,n)    
& 1.6  & 3.1  & 7.7  & 22.8  & 79   
 \\ \hline     
  R_4(B,n)    
& 10  & 28  & 46  & 100  & 231   
 \\ \hline     
  R/E    
& 6.32  & 8.95  & 6  & 4.39  & 2.92  
 \\ \hline \end{array} $$    
 \tabcaption{ Comparaisons du  mod\`ele heuristique avec les valeurs
  exp\'erimentales pour  $d=4$, $n=40$}

\bul $n=45$, $C_4(45)=3.242\cdot 10^{-7}$,

$$  \begin{array}{|c|c|c|c|c|c|} \hline k     
& 1  & 2  & 4  & 8  & 16   
 \\ \hline     
  N_4(B,n)    
&3.5\cdot10^{12}  &5\cdot10^{13}  &8.1\cdot10^{14}  &1.5\cdot10^{16}   &2.9\cdot10^{17}  
  \\ \hline  
  P_4(B,n)    
&1.5\cdot10^{-14}  &1.7\cdot10^{-15}  &2.2\cdot10^{-16}  &3.1\cdot10^{-17}   &4.7\cdot10^{-18}   
 \\ \hline    
  E_4(B,n)    
& 0.1  & 0.1  & 0.2  & 0.4  & 1.3   
 \\ \hline     
  R_4(B,n)    
& 1  & 1  & 2  & 4  & 11  
 \\ \hline     
  R/E    
& 20  & 11.9  & 11.5  & 9.11  &  8.46 
  \\ \hline \end{array} $$
\tabcaption{ Comparaisons du  mod\`ele heuristique avec les valeurs
  exp\'erimentales pour  $d=4$, $n=45$}

{\bf d = 5} :

\bul $n=10$, $(C_5(10)=8.250\cdot 10^{-3})$,
$$  \begin{array}{|c|c|c|c|} \hline k    
& 2  & 4  & 8   
 \\ \hline     
  N_5(B,n)    
&4.5\cdot10^{13}  &1.4\cdot10^{15}  &4.6\cdot10^{16}   
 \\ \hline     
  P_5(B,n)    
&3.9\cdot10^{-4}  &9.4\cdot10^{-5}  &4.1\cdot10^{-5}   
 \\ \hline     
  E_5(B,n)   
& 1.8\cdot 10^{10}  & 1.3\cdot 10^{11}  & 1.9\cdot 10^{12}  
 \\ \hline   
  R_5(B,n)    
& 3.0 \cdot 10^{9}   & 5.3 \cdot 10^{10}  & 1.0 \cdot 10^{12}
 \\ \hline 
  R/E    
& 0.17  & 0.41  & 0.55   
 \\ \hline \end{array} 
$$    
\tabcaption{ Comparaisons du  mod\`ele heuristique avec les valeurs
  exp\'erimentales pour  $d=5$, $n=10$}

\bul $n=20$, $(C_5(20)=1.301 \cdot 10^{-4})$,

$$  \begin{array}{|c|c|c|c|} \hline k     
& 1  & 2  & 4   
 \\ \hline     
  N_5(B,n)    
&2.2\cdot10^{15}  &6.5\cdot10^{16}  &2\cdot10^{18}   
 \\ \hline     
  P_5(B,n)    
&5.2\cdot 10^{-9}  &1.8\cdot 10^{-9}  &6.7\cdot 10^{-10}   
 \\ \hline     
  E_5(B,n)    
& 1.1\cdot 10^{ 7}  & 1.2\cdot 10^{8}  & 1.3\cdot 10^{9} 
 \\ \hline     
  R_5(B,n)    
& 24\,937\,474  & 294\,068\,263  &   3\, 575\, 157\, 628
 \\ \hline     
  R/E    
& 2.19  & 2.52  & 2.75 
 \\ \hline \end{array} $$ 
\tabcaption{ Comparaisons du  mod\`ele heuristique avec les valeurs
  exp\'erimentales pour  $d=5$, $n=20$}

	Pour commenter  la comparaison entre les formules heuristiques
        et les r\'esultats exp\'erimentaux ci-dessus, il faut analyser
        la structure (multiplicative) des formules heuristiques : des
        constantes et des facteurs avec un exposant qui peut
        lui-m\^eme \^etre s\'epar\'e entre son terme principal et des
        termes secondaires.

	Pour un premier commentaire, on notera que l'on multiplie un "tr\`es grand" nombre $N_d( B,n)$ par un "tr\`es petit" nombre $P_d( B,n)$ : si le premier est convenablement connu, le second est obtenu apr\`es un certain nombre d'approximations et de simplifications heuristiques. Il se manifeste alors ce que l'on pourrait appeler un ``effet'' exposant : la moindre variation sur l'heuristique du calcul d'un exposant entra\^{\i}ne des variations (multiplicatives) consid\'erables sur le r\'esultat final $E_d(B,n)$.
	
	On constate entre les valeurs heuristiques et les valeurs exp\'erimentales une discordance d'un facteur de l'ordre de 10 --- en fait entre 3 et 6 dans les zones les plus centrales de la confrontation ---, ce qui n'est pas \'enorme compte tenu des facteurs du produit (jusqu'\`a $10^{17}$ pour le "tr\`es grand" nombre, et jusqu'\`a $10^{-17}$ pour le "tr\`es petit" ). Mais surtout, il y a une tr\`es grande  stabilit\'e lorsque $B$ et $n$ varient, ce qui nous permet d'affirmer que la "forme" des formules est valid\'ee, et que les termes principaux des exposants sont \'egalement valid\'es.
	
		Il reste \`a l'\'evidence des am\'eliorations \`a apporter aux constantes (multiplicatives) et aux termes secondaires des exposants, notamment pour supprimer la variation "convexe" du quotient $R/E$ en fonction de $k$ (et donc de $B$). Mais nous ne pensons pas qu'il y ait aujourd'hui un enjeu suffisant pour se lancer dans ce travail d'am\'elioration, d\'elicat (et hasardeux, si l'on ose dire).
	On peut aussi noter que notre \'evaluation heuristique de la
        probabilit\'e $P_d(B,n)$ a \'et\'e guid\'ee par le
        comportement des polyn\^omes dans la "zone centrale" des
        $\F_B$  mais qu'elle pourrait avoir \'et\'e surestim\'ee \`a
        la p\'eriph\'erie (qui compte beaucoup en volume). Cela
        pourrait expliquer, ou expliquer en partie, l'\'ecart\ldots De
        fait, nous constatons exp\'erimentalement cette surestimation,
        et si la raison est bien celle qui vient d'\^etre avanc\'ee,
        il suffirait de r\'eduire les parties $\F_B$ (sous r\'eserve
        que $B$ reste la borne g\'en\'eriquement atteinte) pour mieux
        cibler la zone centrale et am\'eliorer l'accord entre
        heuristique et exp\'erimentation num\'erique. Mais l\`a
        encore, nous ne pensons pas qu'il y ait aujourd'hui un enjeu
        suffisant.

\subsection{ Le mod\`ele dans le cas des polyn\^omes \`a coefficients rationnels}

L'\'etude effectu\'ee dans le cas des polyn\^omes \`a coefficients entiers peut \^etre facilement adapt\'ee au cas de polyn\^omes \`a coefficients rationnels. Il est sans int\'er\^et de la r\'e\'ecrire au complet et nous contenterons d'exposer bri\`evement les trois points qui sont \`a remanier, et de montrer comment la conclusion est l\'eg\`erement modifi\'ee. 

\bigskip 

{\bf Facteur arithm\'etique}

Il est ais\'e de voir que la fonction $S_k(x)=x(x-1)\ldots(x-k+1)/k!$,
pour $x\in\Z$ est p\'eriodique modulo $p$ de p\'eriode
$p^{t+1}$ o\`u $p^t$ est la plus grande puissance  de $p$ qui divise $k!$ (autrement dit $p^t|| k!$).\\
Par cons\'equent la fonction 
$$P(x)=b_0 S_d(x)+b_1S_{d-1}(x)+\cdots+b_{d-1}S_1(x)+b_d,$$
est p\'eriodique modulo $p$ de p\'eriode $p^{t+1}$ o\`u  $p^t|| d!$.
\\
Notons $D'_d(p)$ la probabilit\'e pour un polyn\^ome g\'en\'erique $P$
du pav\'e $\F_B$ de ne pas avoir $p$ comme d.p.p. On a alors \\
- pour $p>d$, $t=0$  et un calcul analogue au cas  des
polyn\^omes \`a coefficients entiers conduit \`a $D'_d(p)=D_d(p)$,\\
- pour $p\le d$, $D'_d(p)\le D_d(p)$ et sa valeur est calcul\'ee num\'eriquement.

Dans l'\'evaluation de l'expression   $P_d(B,n,m)$, il suffira de
remplacer $D_d(p)$ par $D'_d(p)$.

Il faudra alors remplacer $C_d(n)=\prod_{p\le n} D_d(p)$ par
$C'_d(n)=\prod_{p\le n} D'_d(p)$.

\bigskip 

{\bf Volume du pav\'e $\mathcal{F}_B$}
\\[2ex]
On rappelle que les polyn\^omes $P$ des pav\'es
$\F_B$  v\'erifient entre autres deux propri\'et\'es : sur l'intervalle $[-n,n]$ on a $|P(x)| \le B$, et le terme constant $b_d$ est inf\'erieur \`a $B^{3/4}$ en valeur absolue.
\\
Pour $$P(x)=b_0 S_d(x)+b_1S_{d-1}(x)+\cdots+b_{d-1}S_1(x)+b_d,$$
on a par la normalisation
 	$1 \le b_0 \le \frac{B}{S_d(n)}=\frac{B}{\binom{n}{d}}$, 
	$0 \le b_1 \le \frac{1}2 b_0$ et  
     $|b_j| \le \frac{B^{1-\frac{j}{4d}}}{S_{d-j}(n)}=\frac{B^{1-\frac{j}{4d}}}{\binom{n}{d-j}}$.
\\
Cela donne un nombre de polyn\^omes rationnels $N'_d(B,n)$ \\
$$
N'_d(B,n)\approx \frac{B^2}{8\binom{n}{d}^2} \times \frac{2B^{1-\frac2{4d}}}{\binom{n}{d-2}}\times\ldots
\times \frac{2B^{1-\frac{d}{4d}}}{\binom{n}{d-d}},
$$
soit encore
 $N'_d(B,n)\approx 2^{d-4}
 B^{\frac{7d}8+\frac78+\frac1{4d}}\left(\binom{n}{d}^2\binom{n}{d-2}\ldots\binom{n}{1}\right)^{-1}$.
\\[1ex]
\bigskip

{\bf Estimation de l'esp\'erance $E'_d(B,n)$ et conclusion}
\\[1ex]
Dans la formule estimant l'esp\'erance, il faut remplacer $C_d(n)$ par $C'_d(n)$ et 
$N_d(B,n)$ par $N'_d(B,n)$. Globalement l'esp\'erance
$E'_d(B,n)=P'_d(B,n)N'_d(B,n)$ pour un m\^eme triplet $(d,B,n)$ va
augmenter par rapport au cas des polyn\^omes \`a coefficients entiers car l'augmentation due \`a $N'_d(B,n)$ est plus forte que la diminution due \`a 
$C'_d(n)$, ce qui va am\'eliorer  les perfomances \`a $d,B,n$ fix\'es.
La comparaison num\'erique avec le cas des polyn\^omes \`a
coefficients entiers est illustr\'ee par le tableau suivant, o\`u
pour $n=50$ et $d=4,5,6$, on donne successivement le rapport $C'_d(n)/C_d(n)$, 
le rapport des estimations pour  $N_d(B,n)$ et
$N'_d(B,n)$ (ind\'ependant  de $B$) et enfin le rapport
$E'_d/E_d$ qui en r\'esulte.

\renewcommand{\arraystretch}{1.6}
$$
\begin{array}{|c|c|c|c|}
\hline
  d & 4 & 5 & 6 \\
\hline
C'_d/C_d  &  1.4\cdot 10^{-2} &  4.6\cdot10^{-4}  & 1.2\cdot10^{-4}\\
\hline
N'_d/N_d & 1.9\cdot10^{2} & 2.8\cdot10^{4} &  2.8\cdot 10^{7} \\
\hline
E'_d/E_d&  2.6 & 1.3 \cdot10 & 3.4 \cdot10^{3}\\
\hline
\end{array}
$$
 \tabcaption{Comparaison de l'estimation de l'esp\'erance dans le cas des
   poyn\^omes \`a coefficients entiers ou rationnels}
     
\subsection{Heuristique des pr\'evisions.} 

\subsubsection*{Cas des coefficients entiers}   

	Si l'on veut utiliser notre mod\`ele heuristique pour effectuer des pr\'evisions, il faut adopter la m\^eme probl\'ematique que celle qui a conduit \`a la d\'efinition du mur de Schinzel : pour $d$ et $n$ donn\'es, quelle est la valeur de $B= B_c\ge n^d$  telle que $E_d(n, B_c ) = 1$ ?

	C'est un simple probl\`eme de calcul, un peu long car les formules "compl\`etes" sont compliqu\'ees.

	La recherche exp\'erimentale sugg\`ere que la valeur critique de $B$ est de l'ordre de $n^d$ ou un peu sup\'erieure.
 Si on fixe $d$ et que l'on prend $B=n^d$, lorsque l'on fait cro\^{\i}tre $n$, il existe une valeur critique $n=n_c(d)$ \`a partir de laquelle on a 
$E_d(n_c(d), n^d) \le1$ ; le tableau ci-dessous donne cette valeur pour les valeurs de $d$ qui correspondent aux recherches exp\'erimentales entreprises : 
$$
\begin{array}{|c|c|c|c|c|}
  \hline
d
& 3
& 4
& 5
& 6
\\
\hline
n_c(d)
&  37
& 41
& 47
& 61
\\
\hline
\end{array}
$$
 \tabcaption{Valeurs de la valeur critique $n_c(d)$}

Compte tenu de la forte croissance de $E_d(n, B)$ (facteur principal $B^{7d/8}$), la situation – heuristique – est la suivante :
\begin{itemize}
\item[-] si $n = n_c(d)$, la valeur critique $B_c$  est \'egale \`a $n^d$ ;
\item[-] si $n > n_c(d)$, la valeur critique $B_c$  est l\'eg\`erement sup\'erieure \`a $n^d$.
\end{itemize}

Pour pr\'eciser l'adjectif "l\'eg\`erement", nous allons poser
        $B= n^{td}$ (ce n'est pas la convention que nous avions prise
        pour comparer les valeurs heuristiques et les valeurs
        exp\'erimentales, mais celle-ci est plus commode pour les
        calculs de pr\'evisions).

	Imaginons que l'on veuille entreprendre la recherche de polyn\^omes \`a coefficients entiers de performance sup\'erieure ou \'egale \`a 50, quel est le degr\'e optimal o\`u rechercher ?

R\'eponse avec notre mod\`ele heuristique : le degr\'e o\`u la zone
brute \`a explorer (i.e. avant cribles) est la moins \'etendue. Les
deux tableaux ci-dessous donnent, pour $n = 50$, 55 et 60, la valeur $t_c$ de
$t$ (obtenue num\'eriquement par dichotomie) telle que  $B_c= n^{td}$  et le volume $N_d(B_c,n$) de la zone \`a explorer.

\renewcommand{\arraystretch}{1.6} 

\bul $n = 50$

$$
\begin{array}{|c|c|c|c|c|}
\hline
d
& 3
& 4
& 5
& 6
\\
\hline
t_c &  2.02     &    1.28   &     1.04 &       1.00      
\\
\hline
B_c &
  2.1\cdot 10^{10}   &  4.6\cdot10^{8 } &   6.5\cdot10^{8 }   & 1.6\cdot 10^{10}   
\\
\hline
N_d(B_c,n) &
 3.7\cdot10^{25}  &  4.6\cdot10^{20}  &  1.2\cdot 10^{21}  &  3.9\cdot10^{27} 
\\
\hline
\end{array}
$$
\tabcaption{Pr\'evisions des valeurs critiques  pour $n=50$}

La r\'eponse est donc a priori  pour $n=50$ : $d=4$ ou 5.\\
Notons que l'on a trouv\'e une performance de 46 sur 46 en degr\'e 4 
et une performance de 49 sur 49 en degr\'e 5.

{\bf Une id\'ee du temps de calcul n\'ecessaire.}  Avec les performances actuelles des ordinateurs (juin 2011) , on obtient par exemple en degr\'e 5 ou 6, un temps de calcul de 12 s environ pour explorer $10^{12}$ polyn\^omes (compt\'es avant cribles). Cela donne  pour un seul processeur 140 jours pour explorer une zone de $10^{18}$ polyn\^omes et 38 ans (!) pour une zone de $10^{20}$ polyn\^omes.

\bul $n = 55$

$$
\begin{array}{|c|c|c|c|c|}
\hline
d
& 3
& 4
& 5
& 6
\\
\hline
t_c &     2.28  &     1.43  &    1.11   &    1.00   
\\
\hline
B_c &
  7.6\cdot 10^{11}   &  8.3 \cdot 10^{9 } &   4.8\cdot10^{9 }   & 2.8 \cdot10^{10}   
\\
\hline
N_d(B_c,n) &
 7.2\cdot 10^{30}  &  6.0 \cdot 10^{25}  &  6.2\cdot 10^{24}  &  1.6\cdot 10^{28} 
\\
\hline
\end{array}
$$
\tabcaption{Pr\'evisions des valeurs critiques  pour $n=55$}

La r\'eponse pour $n=55$ est a priori   $d=5$  mais le nombre de polyn\^omes \`a explorer rend  vaine cette recherche.

\bul $n = 60$

$$
\begin{array}{|c|c|c|c|c|}
\hline
d
& 3
& 4
& 5
& 6
\\
\hline
t_c &  2.53   &      1.53     &    1.18    &     1.00  
\\
\hline
B_c &
  2.9\cdot 10^{13}    & 7.0\cdot 10^{10}    & 2.8\cdot 10^{10}   &  5.0\cdot 10^{10}   
\\
\hline
N_d(B_c,n) &
  2.0\cdot 10^{36}    &  3.0\cdot 10^{29}  &   1.6\cdot 10^{28}   &   8.3\cdot 10^{28 }   
\\
\hline
\end{array}
$$
\tabcaption{Pr\'evisions des valeurs critiques pour $n=60$}

La r\'eponse est donc a priori  pour $n=60$  encore $d=5$  mais le nombre de polyn\^omes \`a explorer est astronomique.

\subsubsection*{\bf Pr\'evisions  dans le cas des polyn\^omes \`a coefficients rationnels}

On reprend pour $d\ge4$ le m\^eme type de calcul avec les formules modifi\'ees pour le cas rationnel et on obtient les pr\'evisions suivantes. On notera que le coefficient $t$ peut maintenant \^etre inf\'erieur \`a 1 puisque la valeur plancher de $a_0=1$ correspond \`a $B=\binom{n}{d}\le n^d$.

\bul  $n = 50$
$$
\begin{array}{|c|c|c|c|c|}
\hline
d
& 4
& 5
& 6
& 7
\\
\hline
t &        1.26    &      0.99    &      0.82     &    0.74         
\\
\hline
B 
   &  3.4\cdot 10^{8 } &   2.8\cdot 10^{8 }   & 2.3 \cdot 10^{8}   & 6.2 \cdot 10^{8}
\\
\hline
N_d(B,n) 
&  2.4\cdot 10^{22}  &  2.1\cdot 10^{23}  &  2.2\cdot 10^{23} &  1.2\cdot 10^{26} 
\\
\hline
\end{array}
$$
 \tabcaption{Pr\'evisions des valeurs critiques pour $n=50$ dans le
   cas des polyn\^omes \`a coefficients rationnels }

La r\'eponse est donc ici  pour $n=50$,  $d=4$ . 
Notons que l'on a trouv\'e pr\'ecis\'ement deux performances de 49 sur 49 en degr\'e 4, 6 performances d'au moins 50 sur 50 en degr\'e 5 et plus de 50 performances en degr\'e 6.

Il faut temp\'erer les nombres bruts de polyn\^omes \`a tester donn\'es dans le tableau car avec
des coefficients rationnels, on peut facilement \'eviter les petits d.p.p. en
choisissant convenablement d\`es le d\'epart les coefficients
$b_i$. Cela fait chuter consid\'erablement le temps de calcul. Par
exemple en degr\'e 6, le temps de calcul pour un m\^eme nombre de
polyn\^omes est divis\'e par environ
6000 par rapport au cas des polyn\^omes \`a coefficients entiers.
Cela compense l'augmentation du nombre de polyn\^omes \`a tester et
expliquer les "bonnes" performances obtenues.

\bul $n = 55$
$$
\begin{array}{|c|c|c|c|c|}
\hline
d
& 4
& 5
& 6
& 7
\\
\hline
t &          1.38    &      1.07    &      0.87      &    0.77        
\\
\hline
B &
   4.1\cdot 10^{9 } &   1.9\cdot 10^{9 }   & 1.1 \cdot 10^{9}   & 2.6\cdot 10^{9}
\\
\hline
N_d(B,n) &
 4.9\cdot 10^{26}  &  1.1\cdot 10^{27}  &  4.7\cdot 10^{26} &  1.7\cdot 10^{29} 
\\
\hline
\end{array}
$$
 \tabcaption{Pr\'evisions des valeurs critiques pour $n=55$ dans le cas des polyn\^omes \`a coefficients rationnels}
   
La r\'eponse est donc ici  pour $n=55$, $d=4$ ou 6. 
Du c\^ot\'e exp\'erimental,  la premi\`ere performance trouv\'ee d'au moins 55  sur 55 a \'et\'e en degr\'e 5 mais nous n'avons pas explor\'e les quelques  $4.9.10^{26}$ polyn\^omes en degr\'e 4.

\bul  $n = 60$
$$
\begin{array}{|c|c|c|c|c|}
\hline
d
& 4
& 5
& 6
& 7
\\
\hline
t &           1.51    &      1.14    &      0.92      &    0.81        
\\
\hline
B &
   5.7\cdot 10^{10} &   1.3\cdot 10^{10 }   & 5.9 \cdot 10^{9}   & 1.1 \cdot 10^{10}
\\
\hline
N_d(B,n) &
  2.2\cdot 10^{31}  &  7.4\cdot 10^{30}  &  1.6\cdot 10^{30} &  3.2\cdot 10^{32} 
\\
\hline
\end{array}  
$$
  \tabcaption{Pr\'evisions des valeurs critiques pour $n=60$ dans le
    cas des polyn\^omes \`a coefficients rationnels}
   
 La r\'eponse est donc ici  pour $n=60$,   $d=6$. 
Notons que l'on a trouv\'e pr\'ecis\'ement notre record de 58 sur 58 en degr\'e 6.  C'est probablement encore dans ce degr\'e qu'il faudrait chercher le 60 sur 60 si cela \'etait n\'ecessaire.

\section{Et le cas "$k$ sur $n$" ?} 

  	Il semble opportun de rappeler ce que Dress et Olivier avaient not\'e dans \cite{Dre} : "Les performances d'un polyn\^ome apparaissent comme la "somme" d'un terme fonction de sa constante de Hardy-Littlewood et d'un terme d'apparence al\'eatoire (qui est en quelque sorte un pseudo-al\'ea produit par une situation d\'eterministe mais extr\^emement complexe)".

	Pour le cas "$n$ sur $n$", l'apparence al\'eatoire est totale,
        que ce soit en degr\'e 2 ou en degr\'e sup\'erieur. Cela
        explique le succ\`es du mod\`ele heuristique probabiliste qui
        a \'et\'e d\'evelopp\'e dans ce cet article.

Pour le cas "$k$ sur $n$", et si $n$ est "grand", le terme
d\'eterministe est pr\'epond\'erant, mais le pseudo-hasard produit par
la complexit\'e de la situation se traduit par des performances qui
ont l'apparence du ré\'esultat d'une \'epreuve binomiale dont le
param\`etre est proportionnel \`a la constante de Hardy-Littlewood
$C(P)$ du polyn\^ome et inversement proportionnel à $\log |P(x)|$. 
Cette heuristique, pleinement confirm\'ee par les r\'esultats exp\'erimentaux de Dress et Olivier en degr\'e 2,  est confort\'ee (toujours en degr\'e 2) par notre \'etude num\'erique d'un "tr\`es grand" polyn\^ome de Jacobson et Williams (cf. plus haut, § 3.2). Il ne fait aucun doute qu'elle "fonctionne" pour les degr\'es sup\'erieurs. Pour autant, il n'existe actuellement aucun enjeu qui demanderait l'\'elaboration d'un mod\`ele explicite.

\subsection*{Remerciements} 

Les auteurs tiennent \`a remercier tout particuli\`erement les ing\'enieurs du 
Centre de Calcul Intensif des Pays de Loire et du Groupement de Services
Mathrice du CNRS pour la mise \`a disposition des grappes de calculs
ainsi que pour l'aide apport\'ee \`a la parall\'elisation des programmes.

\noindent Fran\c cois DRESS\\
Institut de Math\'ematiques de Bordeaux\\
UMR CNRS 5251\\
Universit\'e Bordeaux I\\
F-33405 TALENCE Cedex 
\\[2ex]
Bernard LANDREAU\\
Laboratoire Angevin de REcherche en MAth\'ematiques\\
UMR CNRS 6093\\
Universit\'e d'Angers\\
2, bd Lavoisier\\
F-49045 ANGERS Cedex 01

\end{document}